\documentclass[a4paper, 12pt]{paper}

\usepackage{graphicx}
\usepackage{amssymb,amsbsy,amsmath,amsfonts,amssymb,amscd}
\usepackage{xcolor}
\usepackage{titlesec}
\titleformat{\subsection}[runin]{\normalfont\large\bfseries}{\thesubsection}{1em}{}[]

\newtheorem{theorem}{Theorem}[section]
\newtheorem{conjecture}[theorem]{Conjecture}
\newtheorem{lemma}[theorem]{Lemma}
\newtheorem{proposition}[theorem]{Proposition}
\newtheorem{corollary}[theorem]{Corollary}
\newtheorem{definition}[theorem]{Definition\rm}

\title{Distinction of representations via Bruhat-Tits buildings of
  $p$-adic groups}

\subtitle{\`A la m\'emoire de mon ami Fran\c cois Court\`es 1970--2016}
\author{Paul  Broussous}

\begin{document}
\maketitle
\tableofcontents

\section*{Introduction}

 Let $G/H$ be a symmetric space over a non-archimedean local field $F$: $G$ is (the group of $F$-points of)
  a reductive group over $F$ and $H\subset G$ is the subgroup of ($F$-rational) points in $G$
fixed by an involution. A local counterpart of the theory of periods of automorphic forms
on ad\`ele groups is the harmonic analysis on the coset space $G/H$. The irreducible
complex representations $\pi$ of $G$ which contribute to hamonic analysis on $G/H$ are those
representations $\pi$ which embed in the induced representation ${\rm Ind}_H^G \, {\mathbb C}$, where
${\mathbb C}$ denotes the trivial character of $H$.  By
Frobenius reciprocity this amounts to asking that the intertwining space 
${\rm Hom}_G (\pi , {\mathbb C} )$ is non zero. Such representations are
called distinguished by $H$. If $\pi$ is distinguished, a non zero linear form $\Lambda \in 
{\rm Hom}_G (\pi , {\mathbb C} )$ is sometimes called a local period for $\pi$ relative to $H$. 
\medskip

 Among symmetric spaces one has the family of Galois symmetric spaces,
 that is quotients of the form ${\mathbb G}(E)/{\mathbb G}(F)$, where
 $E/F$ is a Galois quadratic extension of $p$-adic fields and $\mathbb
 G$ is a reductive group over $F$. By the conjectural local Langlands
 correspondence an irreducible representation $\pi$ of ${\mathbb
G}(E)$ possesses a Galois parameter $\varphi_\pi$. In \cite{Pr2} Dipendra Prasad
 proposes a ``relative local Langlands correspondence'' of conjectural
 nature: he gives a conjectural list of conditions on the parameter
 $\varphi_\pi$ in order that $\pi$ be distinguished by ${\mathbb
   G}(F)$.   
 \medskip

 Among the irreducible representations of $p$-adic reductive groups,
 one is somehow ``universal''; this is  the Steinberg
 representation. Its definition is uniform and it has nice models of
 geometric nature. It is therefore natural to test Prasad's conjecture
 with this particular representation. In fact in the earlier paper
 \cite{Pr}, Prasad gave a conjecture on the Steinberg representation which
 turns out to be a particular case of the previous conjecture. 
 \medskip

  Let ${\mathbb G}(E)/{\mathbb G}(F)$ be a Galois symmetric space and
  assume that  
  $\mathbb G$ is  quasi-split over $F$. In \cite{Pr} Prasad defines a quadratic character $\epsilon$ of 
  ${\mathbb G}(F)$ and makes the following conjecture.
  \medskip

\noindent {\bf Conjecture} (\cite{Pr} Conjecture 3, p. 77). {\it Let ${{\mathbf S}{\mathbf t}}_E$ be the
  Steinberg representation of ${\mathbb G}(E)$.
  \smallskip

  (a) The intertwining space ${\rm Hom}_{{\mathbb G}(F)} ({{\mathbf S}{\mathbf t}}_E ,
  \epsilon )$ is $1$-dimensional.

  (b) If $\chi \not= \epsilon$ is any other character of
  ${\mathbb G}(F)$, then ${\rm Hom}_{{\mathbb G}(F)} ({{\mathbf S}{\mathbf t}}_E ,
  \epsilon ) =0$. }

  \medskip

  In \cite{BC} the author and F. Court\`es gave a proof of Prasad's conjecture
when ${\mathbb G}$ is split over $F$ and $E/F$ is unramified (actually
there were some other conditions on the group $\mathbb G$ and on the
size of the residue field of $F$,  but they were removed later). 
The aim of this expository work is to explain some of the ideas used
in the proof given in  \cite{BC}. 
\medskip

Let $G$ be a reductive group over a $p$-adic field.
The approach of \cite{BC} is based on the  model of the Steinberg
 representation of $G$ given by  the cohomology of its Bruhat-Tits
 building $X_G$. As a topological space, $X_G$ is a locally compact
 space on which $G$ acts properly (mod center). It is a result of
A.  Borel and J.-P. Serre \cite{BS} that as a $G$-module the top cohomology
space with compact support $H_c^{\rm top}(X_G ,{\mathbb C} )$ is an
irreducible smooth representation of $G$ isomorphic to the Steinberg
representation ${{\mathbf S}{\mathbf t}}_G$. From this result it is easy to construct a model of
the Steinberg representation as a subspace of of the space of complex functions on the
set of chambers of $X_G$. Indeed let ${\mathcal H} (X_G )$ be the space of
harmonic functions on chambers of $X_G$, that is complex functions $f$ satisfying
$$
\sum_{C\supset D} f(C) = 0
$$
\noindent for all codimension $1$ simplex $D$ of $X_G$. Then one has
natural isomorphisms of $G$-modules:
$$
{\rm Hom}_{\mathbb C} ({{\mathbf S}{\mathbf t}}_G ,{\mathbb C} )\simeq
{\mathcal H} (X_G )\otimes \epsilon ', \ \ {{\mathbf S}{\mathbf t}}_G \simeq {\mathcal H}
(X_G )^\infty \otimes\epsilon '
$$
\noindent where ${\mathcal H} (X_G )^\infty$ denotes the space of smooth vectors in the
$G$-module ${\mathcal H} (X_G )$; and where $\epsilon '$ is a certain caracter of
$G$. 
\medskip

 In the case of a Galois symmetric space ${\mathbb G}(E)/{\mathbb
   G}(F)$ satisfying the hypothesis of \cite{BC},  a non-zero equivariant
 linear form $\Lambda \in {\rm Hom}_{{\mathbb G}(F)} ({{\mathbf S}{\mathbf t}}_E ,\epsilon
 )$ is given by
 $$
 \Lambda (f) = \sum_{C\subset X_F} f(C), \ f\in {{\mathbf S}{\mathbf t}}_E \simeq {\mathcal H} (X_E
 )^\infty \
 $$
\noindent where the sum is over those chambers of $X_E :=X_{{\mathbb G}(E)}$ which lie in
$X_F := X_{{\mathbb G}(F)}$ (the building $X_F$ embeds in $X_E$ canonically).
\medskip

 Section 4 of this article will be devoted to the proof of the fact
 that the sum above converges,  for all $f\in {\mathcal H} (X_E )^\infty$, to 
define a non-zero linear form. Our approach here will be different
from that of \cite{BC}; it is based on a new ingredient, namely the Poincar\'e
series of an affine Weyl group, that did not appear in \cite{BC}.
\medskip

We also take the opportunity to give an introductory and pedagogical
 treatmeant of the technical bakground of {\S}4. Namely we start with
 a review of the theory of Bruhat-Tits building (section 1), then we
 state the Borel-Serre theorem and give an idea of its proof (section 2). As an
 exercise we give a complete proof in the case of ${\rm
   GL}(2)$. Section 3 is devoted to the Steinberg representation. We
 review its equivalent definitions and its various models. 
 \medskip

 Originally this article was part of a bigger project joint with
 Fran\c cois Court\`es. Unfortunately Fran\c cois passed away in
 septembre 2016 and I resigned myself to writing on my contribution
 only. However I shall say a few words on  Fran\c cois's contribution
 in {\S}4.5.

\medskip

Throughout this article we shall use the following notation. The symbol $F$ will denote a non-archimedean,
non-discrete, locally compact field. We fix a prime number $p$ and
assume that  $F$ is either a finite
extension of the field ${\mathbb Q}_p$ of $p$-adic numbers or a field
${\mathbb F}_q ((X))$
of Laurent series over a finite field ${\mathbb F}_q$ of $q$ elements,
where $q$ is a power of $p$. For an introduction to  such topological
fields, the reader may read chapter I of \cite{W}, or \cite{Go}. We shall say
that $F$ is a {\it $p$-adic field}. To any $p$-adic field $K$, we
attach: its normalized valuation $v_K$~: $K\longrightarrow {\mathbb Z} \cup
\{+\infty\}$ (assumed to be onto), its valuation ring ${\mathfrak o}_K =\{ x\in
K; \ v_K (x)\geqslant 0\}$, its valuation ideal ${\mathfrak p}_K =\{ x\in
K; \ v_K (x)>0\}$ and its residue field ${\mathbb F}_K ={\mathfrak o}_K /{\mathfrak p}_K$, a finite
extension of ${\mathbb F}_p$. The cardinal of ${\mathbb F}_K$ is
denoted by $q_K$.  We fix a quadratic separable extension
$E/F$. Two cases may occur: either ${\mathfrak p}_F  {\mathfrak o}_E =
{\mathfrak p}_E$ (the
extension is {\it unramified}), or ${\mathfrak p}_F  {\mathfrak o}_E =
{\mathfrak p}_E^2$ (the
extension is ramified). We shall work under the following assumption:
\medskip

(A1) {\it When $E/F$ is ramified, the prime number $p$ is not $2$}.
\medskip

In other words, we assume that the extension $E/F$ is tame (cf. \cite{Fr}
{\S}8). 

 We fix a connected reductive algebraic group  $\mathbb G$ defined over $F$. We
 shall always assume:
\medskip

(A2) {\it The reductive group $\mathbb G$ is split over $F$}.
\medskip

 For simplicity sake, we also assume the following, even though our
 results hold without this assumption:
\medskip

(A3) {\it The root system of $\mathbb G$ is irreducible.}
\medskip

 Prasad's conjecture deal with the symmetric space obtained from the
 reductive group ${\mathbb H}={\rm Res}_{E/F}{\mathbb G}$ (restriction of
 scalars). If $\bar F$ denotes an algebraic closure of $F$, we have an
 isomorphism of $\bar F$-algebraic groups: ${\mathbb H}({\bar F})\simeq
 {\mathbb G}({\bar F})\times  {\mathbb G}({\bar F})$. Let $\sigma$ be
 the $F$-rational involution of ${\mathbb H}$ given by $\sigma (g_1
 ,g_2 )=(g_2 ,g_1 )$; we have  ${\mathbb H}^\sigma = {\mathbb G}$
 (fixed point set). Set
 $G_F ={\mathbb G}(F)$ and $G_E = {\mathbb G}(E)$. Then ${\mathbb
   H}(F)=G_E$ and the action of $\sigma$ on ${\mathbb
   H}(F)$ corresponds to the action of the non-trivial element of the
 Galois group ${\rm Gal}(E/F)$ on $G_E$; this action will be also
 denoted by $\sigma$.  So viewed as a group quotient,
 the symmetric space attached to the group  ${\mathbb H}={\rm
   Res}_{E/F}{\mathbb G}$ equipped with the involution $\sigma$ is
 $G_E /G_F$; this is what we called a Galois symmetric space.

\section{The Bruhat-Tits building} 

\subsection{Apartments and simplicial structure}
 For an introduction to the concept of building the reader read the
 monography \cite{AB}. Basic ideas and various applications of this theory are described in
 \cite{Ro1} and \cite{Ro2}. 
\medskip

 To any reductive group $\mathbb H$ defined over a $p$-adic field $K$,
 the Bruhat-Tits theory (\cite{BT}, \cite{BT2}) attaches a ({\it semisimple}, or
 {\it non-enlarged}) building ${\rm BT}({\mathbb H}, K)$
 equipped with an action of ${\mathbb H}(K)$. In the sequel we
 abreviate $H={\mathbb H}(K)$ and $X_H = {\rm BT}({\mathbb H},
 K)$. Moreover to make things simpler we assume $\mathbb H$ is split
 over $F$ and  that if ${\mathbb Z}$
 denotes the connected center of ${\mathbb H}$, the quotient group
 ${\mathbb H}/{\mathbb Z}$ is simple. We denote by $d$ the $F$-rank of
 that quotient. 
 
\medskip

An outline of the construction of the object $X_H$ is given in
\cite{T}. However, in this expositary paper we shall nearly say nothing  of
this construction. 
\medskip

The $H$-set $X_H$ has a rich structure. First it is a metric space on
which $H$ acts via isometries. Endowed with the metric topology, $X_H$
is locally compact; it is compact (indeed reduced to a single point)
 if and only if the topological group
$H/{\mathbb Z}(K)$ is compact, that is if $d=0$. 

  The set $X_H$ is endowed with a collection of {\it apartments} which
  have the structure of a $d$-dimensional affine euclidean space. They
  play the same r\^ole as charts in differential geometry. More
  precisely $X_H$ is obtained by ``gluing'' these apartments in such a
  way that the following properties are satisfied: 
\medskip

 (1) \hskip1cm  $X_H$ is the union of its apartments,
\smallskip

 (2) \hskip1cm $H$ acts transitively on the set of apartments and if $h\in H$,
 for any apartment ${\mathcal A}$ the induced map ${\mathcal A}\longrightarrow
   h.{\mathcal A}$ is an affine isometry, 
\smallskip

 (3) \hskip1cm for two apartments ${\mathcal A}_1$, ${\mathcal A}_2$, there
 exists $h\in H$ such that $h.{\mathcal A}_1 = {\mathcal A}_2$ and $h$
 fixes ${\mathcal A}_1 \cap {\mathcal A}_2$ pointwise. 
\medskip

We fix a maximal $K$-split torus $\mathbb T$ of ${\mathbb H}$ and
write $T={\mathbb T}(K)$. Let $N(T)$ be the normalizer of $T$ in $H$
and $T^0$ be the maximal compact subgroup of $T$.  The groups $W^\circ
=N(T)/T$ and $W^{\rm Aff}=N(T)/T^0$ are respectively  the spherical
and the {\it extended} affine Weyl groups of $H$ relative to $T$. The group $W^{\circ}$
is a finite reflexion group, indeed a
Coxeter group (cf. \cite{AB}{\S}2). The group $W^{\rm aff}$ is a Coxeter
group if and only if $\mathbb H$ is simply connected as a reductive
$F$-group. In general it may be written as a semidirect product $W^{\rm Aff} = \Omega \rtimes
W_0^{\rm Aff}$, where $\Omega$ is an abelian group and $W_0^{\rm
  Aff}$ is a Coxeter group. 
\medskip

  The torus $\mathbb T$ gives rise to an apartment ${\mathcal A}_T$ of $X_H$
  which is stabilized by $N(T)$. Moreover ${\mathcal A}_T$ is naturally the
  geometric realization of a $d$-dimensional simplicial complex acted
  upon by $N(T)$ via simplicial automorphisms. The maximal dimensional
  simplices of ${\mathcal A}_T$ all have the same dimension $d$. They are called
  {\it chambers}. The subgroup $T^0$ acts trivially on ${\mathcal A}$ so that
  ${\mathcal A}$ is equipped with an action of $W^{\rm Aff}$, whence
  a fortiori of $W_0^{\rm Aff}$. The set of
  chambers
 of ${\mathcal A}_T$ is a principal
  homogeneous space under the action of $W_0^{\rm Aff}$ : for any two
  chambers $C_1$, $C_2$ of ${\mathcal A}_T$, there exists a unique element $w$
  of $W_0^{\rm Aff}$ such that $C_2 =wC_1$. 
\medskip

 The simplicial structure of ${\mathcal A}$ extends in a unique way on the
 whole $X_H$ so that $H$ acts on $X_H$ via simplicial automophisms. A
 simplex of $X_H$ of dimension $d-1$ will be called a {\it codimension
   $1$ simplex}. Each codimension $1$ simplex $D$  of $X_H$ is
 contained in two chambers of ${\mathcal B}$, for any apartment
 ${\mathcal B}$ containing $D$, but is contained in $q_K +1$ chambers
   of $X_H$. For instance when  ${\mathbb H}$ is ${\rm GL}(2)$ or
   ${\rm SL}(2)$, the apartments are euclidean lines, the facets are edges and
   the codimension $1$ facets are vertices. In fact $X_H$ is a
   uniform tree of valency $q_K +1$. 
\medskip

  \subsection{Chambers and Iwahori subgroups} So buildings may also be viewed as combinatorial
  objects  obtained by gluing
  chambers together. Moreover together with properties (1), (2), (3),
  we have:
\medskip

(4) \hskip1cm for any two chambers of $X_H$ there exists an apartment containing
them both. 
\medskip

 From this point of view, it is useful to introduce another distance
 on $X_H$ of combinatorial nature. Two chambers $C_1$ and $C_2$ are
 called {\it adjacent} if the intersection $C_1\cup C_2$ is a
 codimension $1$ simplex. A gallery ${\mathcal G}$ in $X_H$ is a
   sequence ${\mathcal G} =(C_1 ,C_2 ,..., C_d )$ of chambers such
   that, for $i=1,...,d-1$, $C_i$ and $C_{i+1}$ are adjacent. The
   {\it length} of $\mathcal G$ is $d-1$. The {\it combinatorial
     distance} $d(C,C')$ between two chambers $C$, $C'$  is the
   length of a minimal gallery   ${\mathcal G} =(C_1 ,C_2 ,..., C_d )$
   connecting $C$ and $C'$ (i.e. such that $C_1 =C$ and $C_d = C'$).
In fact if $C$ and $C'$ lie in ${\mathcal A}_T$, any minimal gallery connecting
them is contained in ${\mathcal A}_T$. Moreover if $C'=wC$, where $w\in
W_0^{\rm Aff}$ (recall that this $w$ is unique), then  $d(C,C' )
=l(w)$, where $l$~: $W_0^{\rm Aff}\longrightarrow {\mathbb Z}_{\geqslant 0}$
is the length function of the
Coxeter group $W_0^{\rm Aff}$ (e.g. see \cite{AB}, Corollary 1.75). 
\medskip

 The Bruhat-Tits theory attaches to any chamber $C$ of $X_H$ a compact
 open subgroup $I_C$ of $H$ : the {\it Iwahori subgroup} of $H$ fixing
 $C$. If $H_C$ denotes the stabilizer of $C$ in
 $H$, then $I_C$ is a normal subgroup of $H_C^0$, the maximal compact
 subgroup of $H_C$. When $\mathbb H$ is simply connected, one has
 $I_C =H_C^0$, but the containment $I_C \subset H_C^0$ is strict in
 general. If $C$ is a chamber of ${\mathcal A}_T$, then since $T^0\subset I_C$,
 the product set $I_C W^{\rm Aff} I_C$ has a meaning and we have the
{\it  Bruhat-Iwahori decomposition} : 
$$
H= I_C W^{\rm Aff} I_C = \bigsqcup_{w\in W^{\rm Aff}} I_C w I_C\ .
$$

The set $H^0 = I_C W^{\rm Aff}_0  I_C$ is a subgroup of $H$. It is equal
to $H$ when $\mathbb H$ is simply connected. The pair $(I_C ,N)$ is a
{\it $B$-$N$ pair} in $H^0$ and,  as a simplicial complex, $X_H$ is the
building of this $B$-$N$ pair (cf. \cite{AB} {\S}6).  
 \bigskip

 Fix an apartment ${\mathcal A}_T$, attached to a maximal split torus $T$,
 containing $C$. As a Coxeter group $W_0^{\rm Aff}$ is generated by a finite set of
involutions $S$. An involution $s\in S$ acts  on the apartment ${\mathcal A}_T$
as the reflection according to the hyperplan
  containing a certain  codimension $1$  subsimplex $D_s$ of   $C$. This
  codimension $1$ simplex $D_s$ has the form $\{ v_0 , v_1 ,...,v_d
  \}\backslash \{ v_s\}$. One says that $v_s$ is a vertex {\it of type
    $s$} of $C$, and that the {\it opposit} simplex $D_s$ has {\it
    type} $s$ as well. 

 More precisely $W_0^{\rm Aff}$ has a presentation of the form
$$
W_0^{\rm Aff}=\langle s\in S\ ; \ s^2 = 1, \ (st)^{m_{st}}=1,  \ s\not= t\in
S\rangle
$$
\noindent where $m_{st}$ is an integer $\geqslant 2$ or is $\infty$ when
$st$ has infinite order. The length function $l$ has the following
interpretation.  If $w\in W_0^{\rm Aff}$, $l(w)$ is the
number  of involutions in any minimal word on the alphabet $S$
representing $w$.

An important feature of buildings is that they are {\it labellable} as
simplicial complexes.  Let $\Delta_d$ be the standard $d$-dimensional
simplex. Its vertex set is $\Delta_d^0 = \{ 0,1,...,d \}$ and any
subset of $\Delta_d^0$ is allowed to be a simplex.  A {\it labelling} of
$X_H$ is a simplicial map $\lambda$~: $X_H\longrightarrow \Delta_d$ which
preserves the dimension of simplices. In other words, the labelling
$\lambda$ attaches a number $\lambda (s)\in \{ 0,1,...,d\}$ (a {\it
  label}) to  any
vertex $s$ of $X_H$, in such a way that if $\{ s_0, ...,s_k \}$ is a
simplex, then the labels $\lambda (s_0 )$, ..., $\lambda (s_k )$ are
pairewise distinct. 

If $\mathbb H$ is simply connected, then the action of $H$ preserves
the labelling. But this is false in general. In any case the action of
$H^0$ is label-preserving. Let $g\in H$ and $C= \{ s_0 ,...,s_d\}$ be
a chambre of $X_H$. We may consider the permutation $\sigma_{g,C}$ in
${\mathfrak S}_{d+1}$ given by 
$$
\sigma_{g,C}=
\left(
\begin{array}{cccc}
\lambda (s_0 ) & \lambda (s_1 ) & \cdots & \lambda (s_d )\\
\lambda (g.s_0 ) & \lambda (g.s_1 ) & \cdots & \lambda (g.s_d )
\end{array}\right)
$$

\noindent Then the signature of $\sigma_{g,C}$ does not depend on the
choice of $C$; we denote it by $\epsilon_H (g)$. The map
$\epsilon_H$~: $H\longrightarrow \{ \pm 1\}$, $g\mapsto \epsilon_H (g)$ is a
quadratic character of $H$. It is trivial when $\mathbb H$ is simply
connected. 
 \bigskip

\subsection{Behaviour under field extensions} Now let $E/F$ a  tame
quadratic
extension of $p$-adic fields and
$\mathbb G$ be a split reductive $F$-algebraic group with irreducible
root system. Write $\sigma$ for the generator of ${\rm Gal}(E/F)$. 
Write $X_F$ for the Bruhat-Tits building of $\mathbb G$
and $X_E$ for the Bruhat-Tits building of $\mathbb G$ considered as an
$E$-group. These are $G_F$-set and $G_E$-set respectively, where we
put $G_F = {\mathbb G}(F)$ and $G_E = {\mathbb G}(E)$. 

We have a natural action of ${\rm Gal}(E/F)$ on $X_E$ (cf. \cite{T}).
In the simply connected case,  the simplest way to construct it
 is as follows.
 Since ${\rm
   Gal}(E/F)$ acts  continuously on $G_E$ it acts on the set of
 maximal compact subgroups of $G_E$. If $s$ is a vertex of $X_E$,
 there is a unique maximal open subgroup $K_s$ of $G_E$ which fixes
 $s$. One defines $\sigma. s$ to be the unique vertex of $X_E$ fixed by
 $\sigma (K_s )$. Then the action of $\sigma$ on the vertex set of
 $X_E$ extends in an unique way to an affine action of $\sigma$ on the
 whole $X_E$ : if $x\in X_E$ lies in a chamber $C= \{ s_0 , s_1
 ,...,s_d\}$ of $X_E$, with barycentric coordinates $(p_0 ,p_1
 ,...,p_d )$, one defines $\sigma . x$ to be the barycenter of the
 weighted system of points $\{ (\sigma .s_0 , p_0 ),..., (\sigma .s_d
 ,p_d )\}$. 

The action of $\sigma$ on $X_E$ is affine, isometric and
 simplicial. Moreover $\sigma$ permutes the apartments of $X_E$. The
 fixed point set $X_E^{{\rm Gal}(E/F)}$ canonically identifies with
 $X_F$ as a $G_F$-set. So we may view $X_F$ as contained in
 $X_E$. This is a convex subset and we may normalize the metrics in
 such a way that $X_F$ is a submetric space of $X_E$. If $\mathbb T$
 is a maximal $F$-split torus of ${\mathbb G}$ then it is a maximal
 $E$-split torus of ${\mathbb G}$ considered as an $E$-group. Then the
 associate apartments ${\mathcal A}_{{\mathbb  T}(F)}\subset X_F$ and
 ${\mathcal A}_{{\mathbb T}(E)}\subset X_E$ coincide. In particular $X_F$ and
 $X_E$ have the same dimension. 

 If $E/F$ is unramified, then $X_F$ is a subsimplicial complex of
 $X_E$. However if $E/F$ is ramified, the inclusion $X_F\subset X_E$
 is not simplicial. In fact in this case, if $d$ is the dimension of $X_F$, any
 chamber of $X_F$ is the union of $2^d$ chambers of $X_E$. 

 If the extension $E/F$ is not tame, then one still has an  embedding
 $X_F\subset X_E$ which is $G_F$-equivariant, affine and
 isometric. The subset $X_F$ lies in the set  ${\rm Gal}(E/F)$-fixed
 points in $X_E$, but this latter set is strictly larger.

 \subsection{The building of ${\rm GL}(n)$}
 We now work out the example of ${\mathbb H}={\rm GL}(n)$, where
  $n\geqslant 2$ is a fixed integer (references for  more reading are
  \cite{AB}{\S}6.9 and \cite{Ga}{\S}{\S}18, 19). Here the group of $K$-points of
  the  connected center is $Z\simeq K^\times$ and the building
  $X_H$ is of dimension $d=n-1$. In fact the groups ${\rm GL}(n)$,
  ${\rm PGL}(n)$ and ${\rm SL}(n)$ have the same semisimple building. 
\smallskip

  To describe the structure of $H= {\rm GL}(n,K)$ and of its building,
  one makes it act on the $K$-vector space $V= K^n$.  We describe
  first the spherical and affine Weyl groups. We denote by
  $(e_1 ,...,e_n )$ the standard basis of $V$.  As the group of rational
  points of a maximal $K$-split torus, one takes the diagonal torus $T$
  formed of those elements in $G$ that stabilize each line $L_i =
  Ke_i$, $i=1,...,n$. Its normalizer $N$ is  the set of elements
  permutings the lines $L_i$, $i=1,...,n$, i.e. the set of {\it
    monomial} matrices\footnote{A matrix is called monomial if each
    row or column exactly contains a single non-zero coefficient.}.
  The spherical Weyl group $W^\circ$ is isomorphic to the symmetric
  group ${\mathfrak S}_n$. In fact ${\mathfrak S}_n$ embeds in ${\rm
    GL}(n,K)$ in the traditional way so that $N(T)$ is the semidirect
  product $T\rtimes {\mathfrak S}_n$. 
\smallskip

 The group $T^0$ is the set of diagonal matrices in ${\rm GL}(n,K)$
 with coefficients in ${\mathfrak o}_K^\times$, the group of units of the ring
 ${\mathfrak o}_K$. Let $D$ denote the group of diagonal matrices whose
 diagonal coefficients are powers of $\varpi_K$. Then the containment
 $D\rtimes {\mathfrak S}_n\subset N(T)$ induces an isomorphism of
 groups $D\rtimes {\mathfrak S}_n \simeq N(T)/T^0 = W^{\rm Aff}$. 

 For $i=1,...,n-1$, let $s_i$ be the element of ${\rm GL}(n,K)$
 corresponding to the transposition $(i\ i+1)\in {\mathfrak S}_n$. Fix
 a uniformizer $\varpi_K$ of $K$ and write
$$
\Pi = \left(\begin{array}{ccccccc}
0 & 1 & 0 & & \cdots & \cdots & 0\\
0 & 0 & 1 & 0 & \cdots & \cdots & 0\\
\vdots & \vdots & \ddots & \ddots & \ddots & & \vdots\\
\vdots & \vdots & & \ddots & \ddots  & \ddots & \vdots\\
      &         & &        &        &         & 0 \\ 
0 & 0 &  & \cdots & \cdots & 0 & 1\\
\varpi_K & 0 &   & \cdots & \cdots & & 0 
\end{array}\right)
$$
\noindent We put $s_0 = \Pi s_1 \Pi^{-1}$. Then a decompostion $W^{\rm Aff}
= \Omega \rtimes W_0^{\rm Aff}$ is given by $\Omega = \langle \Pi
\rangle$, the group generated by $\Pi$,  and $W_0^{\rm Aff}= \langle
s_0 , s_1 ,...,s_{n-1}\rangle$, the group generated by the $s_i$,
$i=0,...,n-1$ (or more precisely the canonical images of these elements in
$N(T)/T^0$). The $s_i$ are involutions and the group $W_0^{\rm Aff}$
together with the special subset $S= \{ s_0 , s_1 ,...,s_{n-1}\}$ of
generators is a Coxeter system. More precisely, we have the
presentations:
\smallskip

$W_0^{\rm Aff} = \langle s_0 , \ s_1 / \ s_0^2 = s_1^2 =
1\rangle$, if $n=2$,
\smallskip

$W_0^{\rm Aff} = \langle s_0 , \ s_1 , ..., \ s_{n-1}  /
\ s_0^2 = \cdots =s_{n-1}^2 =1, \ (s_i s_{i+1})^3 =1, \ i=0, ..., n-1
\rangle$, if $n\geqslant 3$. Here we have the convention that $s_n
=s_0$.

 Let us now describe the building $X_n$  of ${\rm GL}(n,K)$.  
  A lattice in the $K$-vector space $V=K^n$ is a ${\mathfrak o}_K$-submodule of
  the form $L = {\mathfrak o}_K v_1 + {\mathfrak o}_K v_2 +\cdots +{\mathfrak o}_K v_n$, where $(v_1
  ,v_2 ,...,v_n )$ is a $K$-basis of $V$. Two lattices $L_1$ and $L_2$
  are said equivalent (or homothetic) if there exists $\lambda\in
  K^\times$ such that $L_2 =\lambda L_1$. The equivalence class of a
  lattice $L$ will be denoted by $[L]$. We define a simplicial complex
  ${\mathfrak X}_n$ as follows. Its vertex set is the set of
  equivalence classes of lattices in $V$. A collection of $q+1$
  lattices $[L_0 ]$, $[L_1 ]$, ..., $[L_q ]$ defines a $q$-simplex of
  ${\mathfrak X}_n$ if one can choose the representatives so that 
$$
L_ 0 \supsetneq L_1 \supsetneq L_2 \supsetneq \cdots \supsetneq
L_{q-1}\supsetneq {\mathfrak p}_K L_0\ .
$$

\noindent Then ${\mathfrak X}$ is obviously equipped with an action of
${\rm GL}(n,K)$ via simplicial automorphisms. One can prove \cite{BT2} that the
building $X_n$, as a ${\rm GL}(n,K)$-set, naturally identifies with
the geometric realization of ${\mathfrak X}_n$. 
\smallskip

 In this identification, the vertices belonging to the standard
 apartment ${\mathcal A}_T$ correspond  to the classes $[L]$, where $L$ is a
 lattice {\it split} by the canonical basis of $V$, that is
 satisfying:
$$
L = \sum_{i=1,...,n} L \cap (K\, e_i )
$$
\noindent This identification is compatible with the action of $N(T)$. 
Moreover, if ${\mathcal A}_T^0$ denotes the vertex set of ${\mathcal A}_T$, we have a
surjective map ${\mathbb Z}^n \longrightarrow {\mathcal A}_T^0$, given by 
$$
(m_1 ,...,m_n )\mapsto [\, \sum_{i=1,...,n} {\mathfrak p}^{m_i} e_i \, ]
 \ .
$$
\noindent This map factors through a bijection: ${\mathbb Z}^n /{\mathbb
  Z}\simeq {\mathcal A}_T^0$, where ${\mathbb Z}$ embeds in ${\mathbb Z}^n$ diagonally. As an
euclidean space ${\mathcal A}_T$ is isomorphic to ${\mathbb R}^n /{\mathbb R}$, where ${\mathbb R}$
embeds in ${\mathbb R}^n$ diagonally. The euclidean structure on  ${\mathbb R}^n /{\mathbb R}$
is given as follows: one first equips ${\mathbb R}^n$ with its usual euclidean
structure that one restricts to ${\mathbb R}^n_0 := \{ (x_1 ,...,x_n )\in {\mathbb R}^n
\ / \ x_1 +\cdots + x_n =0\}$; then the quotient  ${\mathbb R}^n /{\mathbb R}$
inherits an euclidean structure via the natural isomorphism of
${\mathbb R}$-vector spaces  ${\mathbb R}^n /{\mathbb R}\simeq {\mathbb R}^n_0$. The action of
$N(T) \simeq T\rtimes {\mathfrak S}_n$ on ${\mathcal A}_T$ is given by 
$$
{\rm diag}(t_1 ,...,t_n ).P_\sigma \ . \ (x_1 , ...,x_n )\ {\rm mod}\ {\mathbb R}
=(x_{\sigma^{-1}(1)} + v_K (t_1 ), ...,
x_{\sigma^{-1}(n)} +v_K (t_n ))\ {\rm mod}\ {\mathbb R}\ ,
 $$
\noindent for all diagonal matrices ${\rm diag}(t_1 ,...,t_n )\in T$ and
all permutation $\sigma \in {\mathfrak S}_n$, where $P_\sigma\in N(T)$
denotes the permutation matrix attached to $\sigma$. 
 The {\it fundamental chamber} in ${\mathcal A}_T$ is the $(n-1)$-simplex $C_0
 = \{
 [L_0 ], ..., [L_{n-1}]\}$, where for $k=0,...,n-1$, $L_k$ is given by 
$$
L_k = \sum_{i=1,...,n-k} {\mathfrak o}_K \, e_i + \sum_{i=n-k+1, ...,n} {\mathfrak p}_K
\, e_i\ .
$$
 
The Iwahori subgroup $I_0$ fixing $C_0$ is called the {\it standard Iwahori 
subgroup} of ${\rm GL}(n,K)$ it is formed of those matrices in ${\rm
  GL}(n,{\mathfrak o}_K )$ which are upper triangular modulo ${\mathfrak p}_K$. The
matrix $\Pi$ stabilizes the chamber $C_0$ : if $i\in \{
0,1,...,n-1\}$, we have $\Pi . [L_i  ]= [L_{i+1}]$, where the index $i$ is
considered {\it modulo} $n$. In fact the stabilizer of $C_0$ in $G$ is
$\langle \Pi \rangle \rtimes I_0$, which is also the normalizer of
$I_0$ in $G$. 
\bigskip

 There is a unique labbeling $\lambda$ on $X_n$ such that $\lambda
 ([L_i ])=i$, $i=0,...,n-1$. It is explicitely given as follows
 (cf. \cite{Ga} {\S}19.3). If $[L]$ is a vertex of $X_n$, choose a
   representative $L$ such that $L\subset L_0$. Since ${\mathfrak o}_K$ is a
   principal ideal domain, the finitely generated torsion
   ${\mathfrak o}_K$-module $L/L_0$ is isomorphic to 
$$
{\mathfrak o}_K /{\mathfrak p}_K^{k_1} \oplus \cdots \oplus  {\mathfrak o}_K /{\mathfrak p}_K^{k_n}
$$

\noindent for some $n$-tupe of integers $0\leqslant  k_1 \leqslant k_2 \leqslant
\cdots \leqslant k_n$. Then 
$$
\lambda ([L])= \sum_{i=0,...,n} k_i\ \ {\rm mod}\ \ n\ .
$$
 
The action of the subgroup $G_0 = I_0 W_0^{\rm Aff} I_0$ of ${\rm
  GL}(n, K)$ preserves the labbelling. In fact the maximal subgroup of ${\rm
  GL}(n, K)$ preserving the labbelling is
$$
\{\  g\in {\rm GL}(n, K), \ \ n\vert v_K ({\rm det}(g))\ \} \ .
$$
The value of the quadratic character $\epsilon = \epsilon_{{\rm
    GL}(n,K)}$ at $\Pi$ is the signature of the cycle $(1\ 2 \ 3
\ \cdots \ n)$, that is $(-1)^{n-1}$. Since ${\rm GL}(n,K)$ is the semidirect
product $\langle \Pi \rangle \rtimes G_0$, it follows that $\epsilon$
is given by
$$
\epsilon (g) = (-1)^{(n-1)v_K ({\rm det}(g))}, \ g\in {\rm GL}(n,K)\ .
$$

\bigskip

 Assume now that $E/F$ is a quadratic  extension of $p$-adic
 fields. Write $X_F$ and $X_E$ for the buildings of ${\rm GL}(n,F)$
 and ${\rm GL}(n,E)$ respectively. The
  containment $X_F\subset X_E$ is given has follows. Set
 $V=F^n$ and identify ${\rm GL}(n,F)$ with ${\rm Aut}_F \, (V)$ and
 ${\rm GL}(n,E)$ with ${\rm Aut}_E (V\otimes_F E)$. Then for any
 ${\mathfrak o}_F$-lattice $L$ of $V$, the vertex $[L]$ of $X_F$ corresponds to
 the vertex $[L\otimes_{{\mathfrak o}_F}{{\mathfrak o}_E}]$ of $X_E$, where
 $L\otimes_{{\mathfrak o}_F}{{\mathfrak o}_E}$ is identified with its canonical image in
 $V\otimes_F E$. 

Let $\mathbb T$ be the diagonal torus of ${\mathbb
   GL}(n)$. It corresponds to an apartment ${\mathcal A}_F$ of $X_F$ and
${\mathcal A}_E$ of $X_E$. We saw that both apartments identify with ${\mathbb R}^n
/{\mathbb R}$. Then the contaiment $X_F \subset X_E$ restricts to ${\mathcal A}_F
\subset {\mathcal A}_E$ (in fact an equality) where it corresponds to the map
$$
{\mathbb R}^n /{\mathbb R} \longrightarrow {\mathbb R}^n /{\mathbb R} , \ \ \  (x_1 ,...,x_n )\ {\rm
  mod}\  {\mathbb R}\mapsto  (e(E/F)x_1 ,...,e(E/F)x_n )\ {\rm mod}\ {\mathbb R} ,
$$
\noindent where $e(E/F)$ is the ramification index of $E/F$.

\vskip1cm
\begin{center}
\includegraphics[scale=0.7]{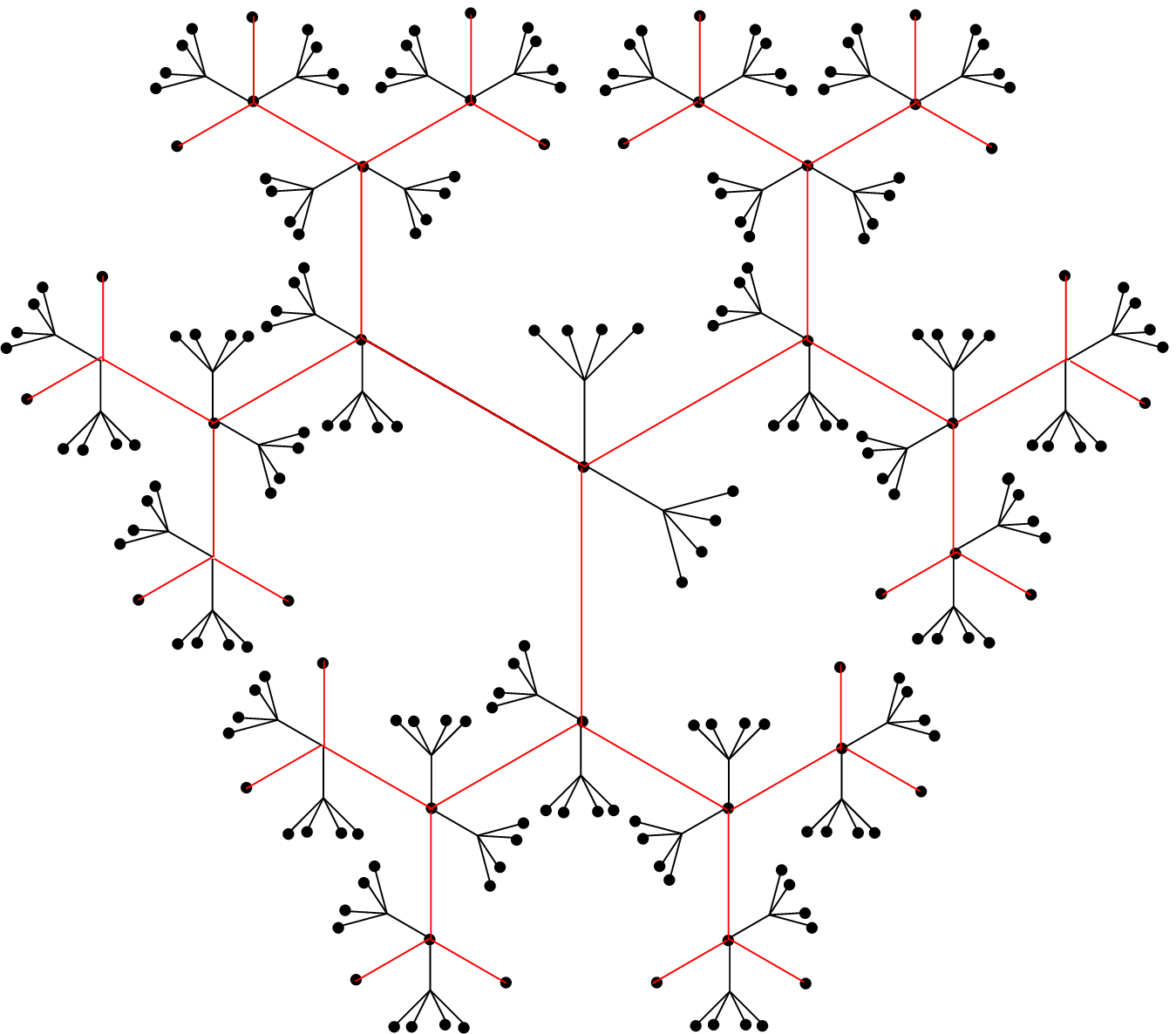}
  \end{center}

\bigskip

\centerline{{\bf Fig. 1.} The embedding $X_F \subset X_E$, $E/F$
  unramified. }
\vskip1cm
\begin{center}
\includegraphics[scale=0.7]{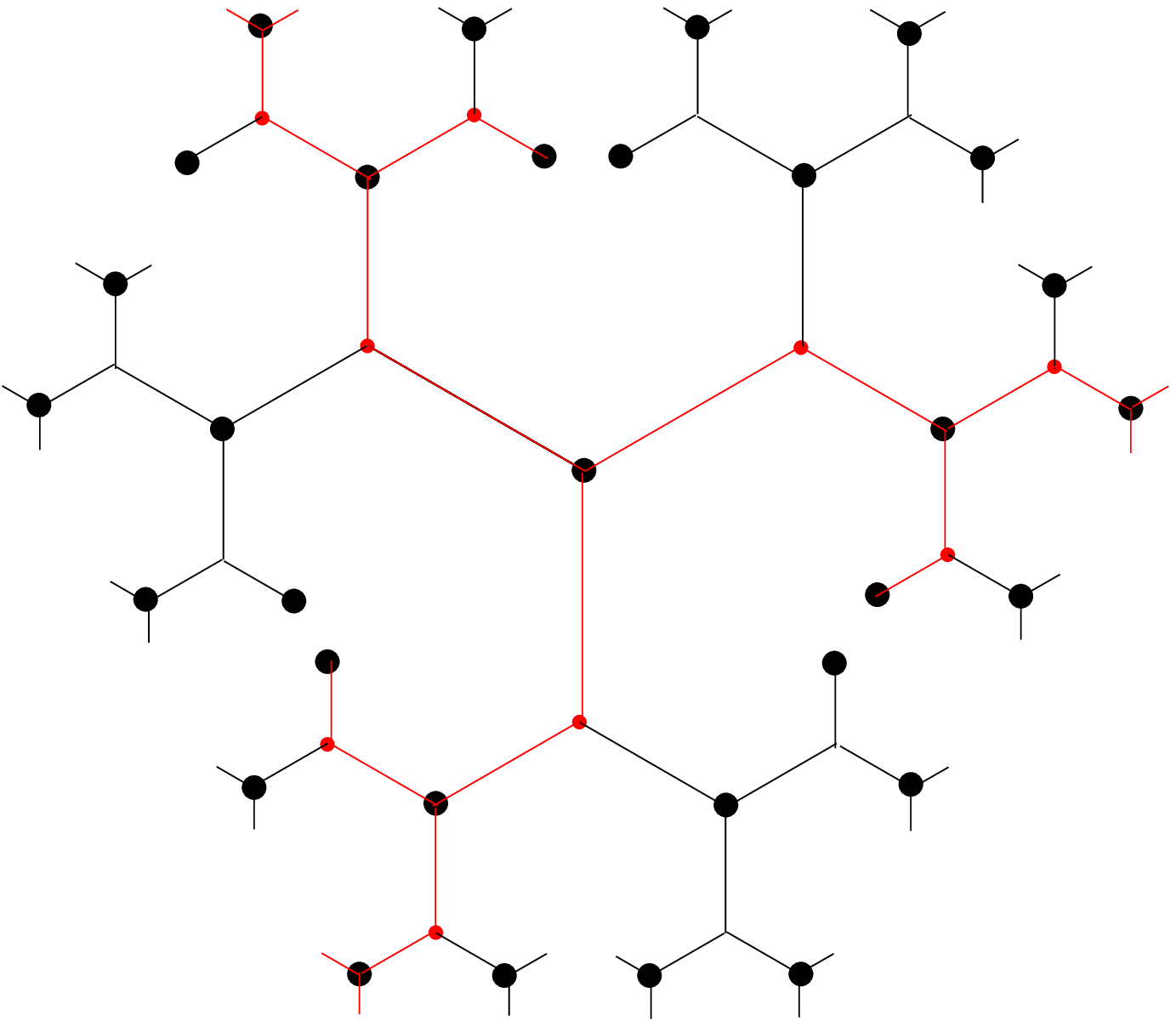}
\end{center}

\bigskip

\centerline{{\bf Fig. 2.}  The embedding $X_F \subset X_E$, $E/F$
  ramified.}

 \vskip1cm

 In figure 1, we drew  part of the building $X_E$ when $F={\mathbb Q}_2$ and
 $E/F$ is unramified. It is an uniform tree of valency $5$. The
 building  $X_F\subset X_E$ is drawn in red; it is an uniform tree of
 valency $3$.  Figure 2 represents a
 part of $X_E$ when   $F={\mathbb Q}_2$ and
 $E/F$ is ramified. Note that in this case $E/F$ is not tame. The
 building of $X_F$ is drawn in red. Both uniform trees $X_F$ and 
$X_E$ have valency $3$. On can see that the embedding
 $X_F\subset X_E$ is not simplicial: the red vertices of $X_E$ are not
 vertixes of $X_F$; they correspond to middles of chambres.

\section{Borel-Serre Theorem}

 \subsection{Statement and ideas of the proof} We fix a $p$-adic field
 $K$ and a split reductive $K$-group $\mathbb
 H$. We use the same notation and assumptions than in {\S}1.1. Recall
 that equipped with it metric topology the affine building $X_H$ is a
 locally compact topological space. So for any integer $k\geqslant 0$,
 we may consider $H_c^k (X_H , {\mathbb C} )$, the cohomology space with compact
 support of $X$ with coefficients in ${\mathbb C}$. Here these cohomology
 spaces are defined by any {\it reasonable} cohomology theory, e.g.
 Alexander-Spanier theory \cite{Ma}, or cohomology in the sense of sheaf
 theory \cite{Bre}. In particular the space $H_c^0 (X_H ,{\mathbb C} )$ is the
 ${\mathbb C}$-vector space of complex locally constant functions with compact
 support on $X_H$. If the dimension $d$ of the building is $>0$, then $X_H$ is
 connected and non-compact, so that $H_c^0 (X_H ,{\mathbb C} )=0$. From now
 on, we assume that $d\geqslant 1$. 
\medskip

 In \cite{BS}, in order to study the cohomology of $S$-arithmetic groups,
 Borel and Serre state and prove the following result. 

\begin{theorem} (Cf. \cite{BS} Th\'eor\`eme 5.6 and {\S}5.10) The cohomology
  space $H_c^k (X_H,{\mathbb C} )$ is trivial when $k\not= d$. Moreover as a
  $H$-module $H_c^d (X_H,{\mathbb C} )$ is smooth and irreducible; it is in fact
  isomorphic to the Steinberg representation of $H$.
\end{theorem} 

 In fact the topological space $X_H$ is contractible. Indeed if $x$,
 $y$ are two points of $X_H$, and $t\in [0,1]$, the barycenter
 $tx+(1-t)y$ is well defined (it is defined in any apartment ${\mathcal A}$
 containing $x$ and $y$ and does not depend on the choice of
 ${\mathcal A}$). Moreover the map
$$
X\times X\times [0,1]\longrightarrow X, \ (x,y,t) \mapsto tx+(1-t)y
$$
\noindent is continuous. So if $o$ is any  point of $X$, the map
$$
F :X\times [0,1] \longrightarrow X, \ (x,t)\mapsto (1-t)x+to
$$
is a homotopy between $F(-,0)$, the identity map of $X_H$, and $F(-,1)$
the constant map with value $o$. 

 It follows that the cohomology space (without support) $H^k (X_H,{\mathbb C} )$
 are trivial when  $k>0$. This also means that if for some $k>0$,
 $H_c^k (X_H ,{\mathbb C} )\not= 0$, this is not due to the existence of ``cycles''
 in $X_H$ but rather to the fact that $X_H$ is not compact. So the
 natural idea that Borel and Serre follow is to compactify the space
 $X_H$ by adding a boundary, and this boundary is the Tits building
 $Y_H$ of $H$ topologized in a certain way that we describe now. 
\medskip

The Tits building $Y_H$ is a simplicial complex. Its 
 vertices are the maximal proper parabolic subgroups of $H$. By
 definition $r+1$ such parabolic subgroups $P_0$, $P_1$, ..., $P_r$
 define a $r$-simplex if the intersection $P_0\cap P_1 \cap \cdots
 \cap P_r$ is a parabolic subgroup of $H$ (or equivalently contains a
 Borel subgroup of $H$). Hence the set of simplices of $Y_H$ is in
 $H$-equivariant bijection with the set of proper parabolic subgroups
 of $H$. The dimension of $Y_H$ is $d-1$, where $d$ is the dimension
 of the affine building $X_H$. Fix a simplex of maximal dimension (a chamber) $D$ of
 $Y_H$. It is a fundamental domain for the action of $H$ on $Y_H$. If
 $B$ denotes the Borel subgroup stabilizing $D$, one may view $Y_H$ as
 a quotient $(H/B) \times D/\sim$  of $(H/B) \times D$.
 On $H/B$ we put the $p$-adic topology
 so that it is a compact set (the group of $K$-points of a complete
 projective variety) and on $D$ we put the euclidean topology (so as a
 simplex, it is compact).  Then the topology on $Y_H$ is the quotient
 topology of $(H/B) \times D/\sim$. As a quotient of a
 compact space $Y_H$ is compact. In particular we have $H^q (Y_H
 ,{\mathbb C}
 )=H^q_c (Y_H ,{\mathbb C} )$, $q\geqslant 0$, where $H^q$ denotes a
 cohomology space without condition of support.   
 \smallskip
 
 The Borel-Serre compactification of $X_H$ is the disjoint union 
${\bar X}_H = X_H \sqcup Y_H$. We shall not describe the topology of
 ${\bar X}_H$. Let us just say that the induced topology on $X_H$
 (resp. $Y_H$) is the metric topology (resp. the topology we defined
 in the last paragraph), that $X_H$ is open and dense in ${\bar X}_H$,
 that $Y_H$ is closed. Moreover as $X_H$, the topological space ${\bar
   X}_H$ is contractible. It follows that its reduced cohomology
 spaces ${\tilde H}_c^q  ({\bar X}_H , {\mathbb C} ) ={\tilde H}^q ({\bar X}_H ,{\mathbb C} )$
 are trivial for all $q$. Recall that ${\tilde H}^q ({\bar X}_H ,{\mathbb C} )$ is
 defined as follows: ${\tilde H}^q ({\bar X}_H ,{\mathbb C} ) = H^q ({\bar X}_H ,{\mathbb C} )$,
 for $q>0$ and ${\tilde H}^0 ({\bar X}_H ,{\mathbb C} ) = H^0 ({\bar X}_H ,{\mathbb C} )/{\mathbb C}$,
 where ${\mathbb C} \subset H^0 ({\bar X}_H ,{\mathbb C} )$ is viewed as the subspace of
 constant functions. 
\smallskip

 In \cite{BS}, the authors prove the following result.

\begin{theorem} The cohomology space  $H^q (Y_H ,{\mathbb C} )$ is trivial for
  $q<d-1$. 
\end{theorem}
 
Moreover they describe ${\tilde H}^{d-1}(Y_H ,{\mathbb C} )$ as a
$H$-module. To state their result we nee to introduce a bit of
notation. We fix a maximal split torus $T$ of $H$,  denote by $\Phi = \Phi
(H,T)$ the corresponding root system. We fix a Borel subgroup $B$
containing $T$  and  a basis $\Delta$ of the set of positive roots in 
$\Phi$ relative to $B$. We have a $1-1$ correspondence $I\mapsto P_I$
between subsets of $\Delta$ and standard parabolic subgroups of $H$
relative to $B$ (in particular $P_\emptyset =B$ and $P_\Delta =
G$). For $I\subset \Delta$ we denote by $\sigma_I$ the representation
of $H$ in  $C^\infty (H/P_I ,{\mathbb C} )$, the space of locally constant
complex functions on the compact $H$-set $H/P_I$. The representations
$\sigma_I$, $I\subset \Delta$, being   smooth and of finite length,
one may consider the following element of the Grothendieck group of
smooth complex representations of $H$ of finite length:
$$
{\rm St}_H := \sum_{I\subset \Delta} (-1)^{\vert I\vert} \, \sigma_I
\ .
$$

This element of the Grothendieck group is actually an irreducible
representation  and is
called the {\it Steinberg representation} (we shall give more details
on this representation in the next section). By exploiting the
combinatorics of the Tits building $Y_H$, Borel and Serre prove that
we have the isomorphism of $H$-modules:
\begin{equation}
{\tilde H}^{d-1}(Y_H ,{\mathbb C} )\simeq {\rm St}_H\ .
\label{CohomTits}
\end{equation}

Now the proof of Theorem 2.1 proceeds as follows. The long exact
sequence of the pair of topological spaces $({\bar X}_H , Y_H )$
writes (\cite{Ma} Theorem 1.6):
\begin{equation}
      \longrightarrow  H_c^{k-1}({\bar X}_H ,{\mathbb C} )  
\longrightarrow H_{c}^{k-1}(Y_H ,{\mathbb C}) \longrightarrow  H_{c}^{k}(X_H ,{\mathbb C} )\longrightarrow  H_{c}^{k}({\bar X}_H ,{\mathbb C} )\longrightarrow
  \cdots , \ k\geqslant 1\ .
\label{longexseq}
\end{equation}
\noindent (The case $k=0$ was already considered above). If $k>1$, then
$H_c^{k-1} ({\bar X}_H  ,{\mathbb C} )$ and $H_c^{k} ({\bar X}_H  ,{\mathbb C} )$
are trivial since ${\bar X}_H$ is contractible. Hence a piece of the
long exact sequence \eqref{longexseq} writes:
\begin{equation}
0\longrightarrow H_{c}^{k-1} ({Y}_H ,{\mathbb C} )\longrightarrow H_c^k (X_H ,{\mathbb C} ) \longrightarrow 0 , 
\end{equation}
\noindent and we obtain the isomorphism of $H$-modules: $ H_c^k (X_H ,{\mathbb C}
)\simeq  H_{c}^{k-1} ({Y}_H ,{\mathbb C} ) =  {\tilde H}_{c}^{k-1} ({Y}_H ,{\mathbb C}
)$, as required. 

\noindent If $k=1$, using again the contractility of ${\bar X}_H$, we
obtain the exact sequence:
$$
{\mathbb C} \simeq H_c^{0}({\bar X}_H,{\mathbb C} ) \overset{j}{\longrightarrow} H_c^0 (Y_H ,{\mathbb C}
) \longrightarrow H_c^1 (X_H ,{\mathbb C} )\longrightarrow 0\ ,
$$
\noindent whence the $H$-isomorphism: $H_c^1 (X_H ,{\mathbb C} )\simeq  H_c^0 (Y_H ,{\mathbb C}
) /j( H_c^{0}({\bar X}_H,{\mathbb C} ))$, where this latter quotient is easily
seen to be isomorphic to  ${\tilde H}_c^0 (Y_H ,{\mathbb C} )$, as required.

\subsection{Sketch of proof for ${\rm GL}(2)$} As an exercice, we give a {\it nearly
  complete}\footnote{In fact easy results or 
  routine details will be left to the reader} and elementary proof of Borel-Serre theorem
for $H={\rm GL}(2,F)$. 

 Recall that $X_H$ is a uniform tree of valency $q_F +1$. The Tits
 building has dimension $1-1=0$, and as a topological space and $H$-set
 it identifies with the quotient $H/B$, where $B$ is the Borel
 subgroup of upper triangular matrices. In turn this quotient
 naturally identifies with the projective line ${\mathbb P}^1 (F)$ (as a
 topological space and $H$-set). Indeed $H$ acts transitively on the
 set of lines of $V=F^2$ and $B$ is the stabilizer of the line
 generated by $(1,0)$. So in this very particular case, Borel-Serre
 theorem claims that we have an isomorphism of $H$-modules: $H_c^1
 (X_H ,{\mathbb C} )\simeq {\tilde H}^0 ({\mathbb P}^1 (F)) =
 C^0 ({\mathbb P}^1 (F))/{\mathbb C}$, where $C^0 ({\mathbb P}^1 (F))/{\mathbb C}$ is
 the space of locally constant complex functions on ${\mathbb P}^1 (F)$
 quotiented by the subspace of constant functions.  Of course the
 $H$-module e $C^0 ({\mathbb P}^1 (F))/{\mathbb C}$ is nothing other than the Steinberg
 representation of $H$.
\bigskip

 In the tree case, the  Borel-Serre compactification coincides with the
 compactification obtained by adding {\it ends} (cf. \cite{Se} I.2.2 and
 II.1.3, \cite{DT} 1.3.4). A half-geodesic in  the tree $X_H$ is a sequence of
 vertices $g=(s_k )_{k\geqslant 0}$ such that for all $k\geqslant 0$,
 $\{ s_k ,s_{k+1}\}$ is an edge ($g$ is a {\it path}) and
 $s_{k+2}\not= s_k$ ($g$  is {\it non-backtracking}). Two
 half-geodesics $b=(s_k )_{k\geqslant 0}$ and $b'=(t_k )_{k\geqslant
   0}$ are said to be {\it equivalent}, if there exists $l\in {\mathbb Z}$,
 such that $s_k =t_{k+l}$, for $k$ large enough. An {\it end} in $X_H$
 is an equivalence class of half-geodesics; we denote by ${\rm End}_H$ this set
 of ends. One observes that $H$ acts naturally on ${\rm End}_H$ and that, a vertex
 $o$ in $X_H$ being fixed, any end has a unique representative $(s_k
 )_{k\geqslant 0}$ such that $s_0 =o$. If $\{ s,t\}$ is an edge of
 $X_H$ we define a subset $\Omega_{(s,t)}$ of ${\rm End}_H$ as follows: an end
 belongs to $\Omega_{(s,t)}$ if  its representative $(s_k )_{k\geqslant
     0}$ with $s_0 = s$ satisfies $s_1 =t$. 
We equip ${\rm End}_H$ with the topology whose basis of open subsets if formed
of the $\Omega_{(s,t)}$, where $(s,t)$ runs over the ordered edges of
$X_H$.   
 \smallskip

 There is a $H$-equivariant homeomorphism $\varphi$~: ${\rm End}_H \overset{\sim}{\longrightarrow}
 {\mathbb P}^1 (F)$ that we describe now. If $b\in {\rm End}_H$, let $g=(s_k
 )_{k\geqslant 0}\in b$ be the representative satisfying $s_0 =[{\mathfrak o}_F
   e_1 +{\mathfrak o}_F e_2 ]$, where $(e_1 ,e_2 )$ is the canonical basis of
 $F^2$. Then (cf. [DT] 1.3.4) one can find a basis $(v_1 ,v_2 )$ of
 $F^2$ such that $s_k =[{\mathfrak o} v_1 +{\mathfrak p}^k v_2 ]$, $k\geqslant 0$. We
 then define $\varphi (b)$ to be the line $Fv_1 \in {\mathbb P}^1
 (F)$. Conversely if $[x \, : \, y] := {\rm Vect}_F \, (x,y)$ is a line in
 ${\mathbb P}^1 (F)$, one may arrange the representatives $x$, $y$ to lie in
 ${\mathfrak o}_F$ and to verify:  $x$ or $y\in {\mathfrak o}_F^\times$. Then the end
 $b$ with representative $(s_k )$ defined by $s_k =[{\mathfrak o}_F (x,y)
   +{\mathfrak p}^k ({\mathfrak o}_F e_1 +{\mathfrak o}_F  e_2 )]$, $k\geqslant 0$, satisfies
 $\varphi (b) = [x\, : \, y]$.  
\smallskip

 In the sequel we canonically identify ${\rm End}_H$ and ${\mathbb P}^1 (F)$
 as $H$-sets and topological spaces. 
\smallskip

\begin{center}
\includegraphics[scale=0.7]{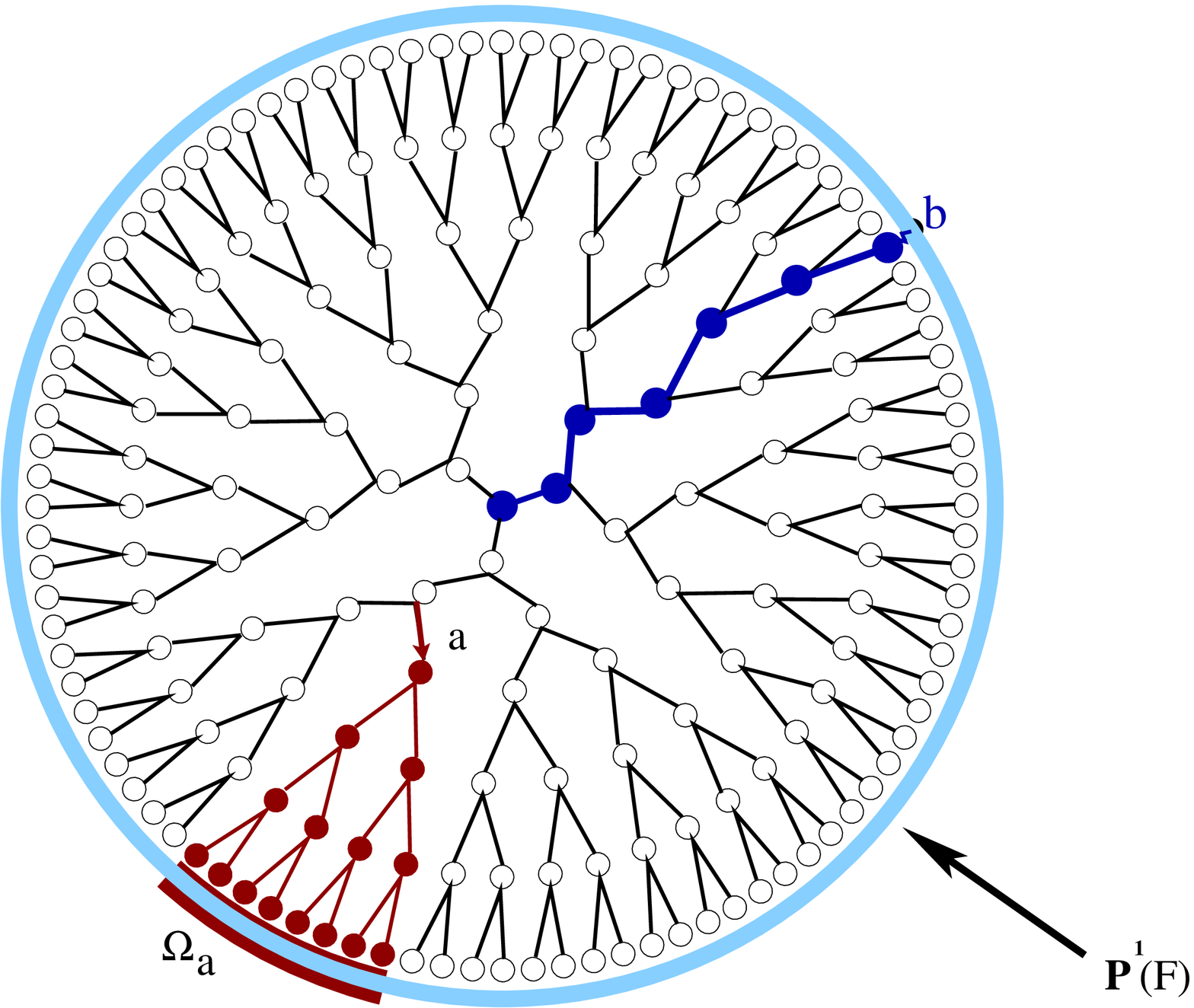}
\end{center}

\centerline{{\bf Fig. 3.}  An end of $X_H$ and a basic open subset of
  ${\mathbb P}^1 (F)$}
\bigskip

 Since $X_H$ is a simplicial complex, its cohomology space $H_c^1 (X_H
 ,{\mathbb C} )$ may be computed via the complex of simplicial cochains. This
 argument will be used again later in these notes.  Write
 $X_H^0$ for the set of vertices of $X_H$, and $X_H^1$ for the set of
 {\it oriented edges} of $X_H$, that is the set of ordered pairs
 $(s,t)$, where $\{ s,t\}$ is an edge of $X_H$. We denote by $C_c^0
 (X_H )$ the ${\mathbb C}$-vector space of {\it $0$-cochains with compact
   support}, that is the set of functions $f$~: $X_H^0 \longrightarrow {\mathbb C}$ with finite
 support. Similarly $C_c^1 (X_H )$ denotes the ${\mathbb C}$-vector space of
 {\it $1$-cochains with compact support}, that is the set of functions
 $\omega$~: $X_H^1 \longrightarrow {\mathbb C}$ with finite support and satisfying
 $\omega (s,t)=-\omega (t,s)$, for all edges $\{ s,t\}$ of $X_H$.  We
 have a coboundary operator $d$~: $C_c^0 (X_H )\longrightarrow C_c^1 (X_H )$ given
 by $df(s,t)=f(s)-f(t)$, for all edges $\{ s,t\}$. The spaces $C_c^i
 (X_H )$ are naturally smooth $H$-modules and the map $d$ is
 $H$-equivariant. As a smooth representation of $H$, the cohomology of
 $X_H$ is given by the cohomology of the complex :
$$
0\longrightarrow C_c^0 (X_H )\overset{d}{\longrightarrow} C_c^1 (X_H )\longrightarrow 0\ .
$$

\noindent In particular we have an isomorphism of $H$-modules : $H_c^1 (X_H
,{\mathbb C} )\simeq C_c^1 (X_H )/dC_c^0 (X_H )$.  
\smallskip

 We are now going going to construct a natural $H$-equivariant map
 $\Psi$~: $C_c^1 (X_H )\longrightarrow C^0 ({\mathbb P}^1 (F))/{\mathbb C}$ and prove that it is
 onto and has kernel $dC_c^0 (X_H )$. The Borel-Serre theorem will
 follow. 
\medskip

 Once for all fix a vertex $o\in X_H^0$ (e.g. $o=[{\mathfrak o}_F e_1 +{\mathfrak o}_F
   e_2 ]$). If $p=(s_0 ,s_1 ,...,s_k )$ is any path in $X_H$ and
 $\omega \in C_c^1 (X_H)$, define the {\it integral of $\omega$ along
   $p$} to be 
$$
\int_p \omega := \sum_{i=0,...,k-1} \omega (s_i ,s_{i+1})\ .
$$

\noindent Note that, since $X_H$ is simply connected, 
 $\int_p \omega$ only depends on $\omega$ and the
origin and end $s_0$ and $s_k$ of $p$. If $\omega\in C_c^1 (X_H )$ and
$b$ is an end of $X_H$ with representative $(s_k )_{k\geqslant 0}$
(normalized by $s_0 =o$), the  sequence $\displaystyle (\int_{(s_0 ,...,s_k
  )}\omega )$ is stationnary since $\omega$ has finite support; denote
by $\varphi_\omega (b)$ its ultimate value. It is a routine exercice to
prove that $f_b$ is a locally constant function on ${\mathbb P}^1 (F)$. 
 We then define $\Psi (b)$ to
be the image of the function $\varphi_\omega$~: ${\mathbb P}^1 (F)\longrightarrow {\mathbb C}$ in $C^0
({\mathbb P}^1 (F))/{\mathbb C}$. Note that if one changes the origin vertex $o$ the
function $\varphi_\omega$ is modified by an additive constant so that its
image in $C^0 ({\mathbb P}^1 (F))/{\mathbb C}$ does not change. 
\smallskip
 
 Let us prove that the kernel of $\Psi$ is $dC_c^0 (X_H )$. The
 containment $dC_c (X_H )\subset {\rm Ker}\, \Psi$ is easy for if
 $f\in C_c^0 (X_H )$, one has $\varphi_{df}(b) =-f(o)$ for any end $b$
 so that $\varphi_{df}$ is constant.  Let $\omega\in {\rm Ker}\,
 \Psi$. This means that $\varphi_\omega$ is a constant function. Let
 $c$ be its value. Then if an end $b$ has representative $(s_k )_{k\geqslant 0}$,
 with $s_0 =o$, one has $\int_{s_0 ,...,s_k }\omega =c$ for $k$ large
 enough. Define a function $f$ on $X_H^0$ by 
$$
f(s) = \left(\int_{p\ : \ o\longrightarrow s} \omega\right)  -c\ , \ s\in X_H^0
$$
\noindent where the notation $p$~: $o\longrightarrow s$ means that $p$ is any path
from $o$ to $s$. Since $X_H$ is a tree, $f$ is well defined, and by a
compactness argument its has finite support. It is finally clear that
$df=\omega$, as required. 
\smallskip
  
 For the surjectivity of $\Psi$, fix $g\in C^0 ({\mathbb P}^1 (F))$ be any
 locally constant function. On has to find $\omega\in C_c^1 (X_H )$
 satisfying $\varphi_\omega =g$. For any integer $r\geqslant 1$, 
consider the finite subtree $S(o,r)$ of formed of  points at distance
less than or equal to $r$ from a fixed vertex $o$.\footnote{The
  distance on the tree is normalized so that the length of an edge is
  $1$.} Then $S(o,r)$ contains 
$$
1+(q_F +1)+(q_F +1)q_F +\cdots + (q_F +1)q_F^i +\cdots + (q_F
+1)q_F^{r-1} = 1+ (q_F +1)\frac{q_F^r -1}{q_F -1}
$$
\noindent  vertices. A vertex of $S(o,r)$ has valency $1$ or $q_F +1$
according to whether it is an {\it end} of $S(o,r)$ or not. Let $s_i$,
$i=1,...,(q_F +1)q_F^{r-1}$ be an indexing of the ends of $S(o,r)$, and
for each $s_i$, let $t_i$ denote the unique neighbour vertex of $s_i$
in $S(o,r)$. Then we have the following partition of ${\mathbb P}^1 (F)$:
$$
{\mathbb P}^1 (F) =\bigsqcup_{i =1,...,(q_F +1)q_F^{r-1}} \Omega_{(t_i ,s_i )}
\ , 
$$
\noindent and by a compactness argument we may assume, by taking $r$ large
enough,  that $g$ is constant
on each $\Omega_{(t_i ,s_i )}$; write $c_i$ for this constant value.  Now define a function $f_r$ on the
set of vertices of $S(o,r)$ by $f_r (s_i ) = c_i$, $i=1,...,(q_F
+1)q_F^{r-1}$, and by giving arbitrary values to $f_r (s)$, for all
vertices of $S(o,r)$ which are not ends. Define $\omega \in C_c^1 (X_F
,{\mathbb C} )$ by $\omega (s,t)=0$, if the edge $\{ s,t\}$ doe not lie in
$S(o,r)$, and by $\omega (s,t)=f_r (t)-f_r (s)$ otherwise. It is easy
to check that $\Psi_\omega =g$, as required.

\begin{center}
\includegraphics[scale=0.7]{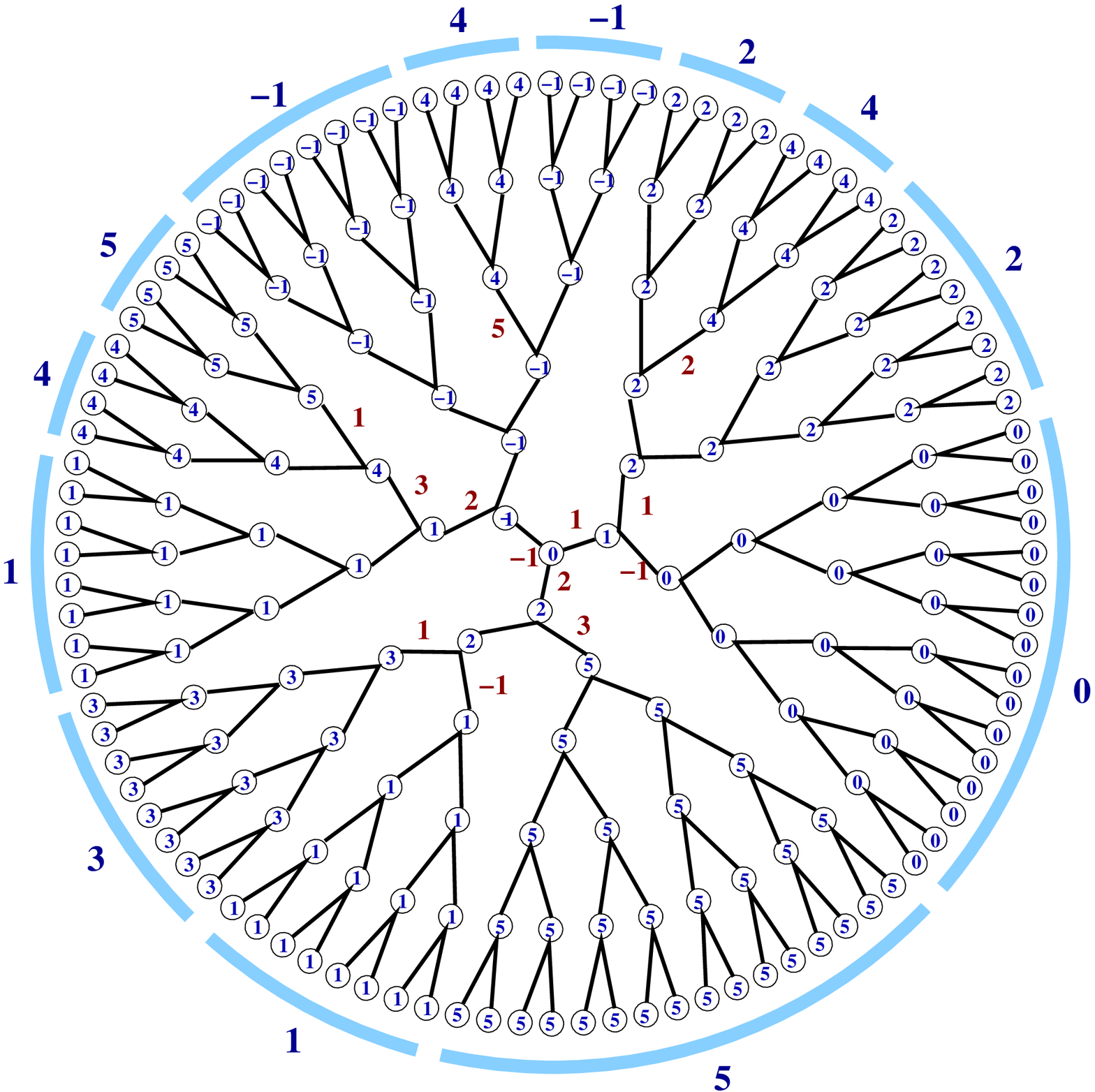}
\end{center}

\centerline{{\bf Fig. 4.}  A $1$-cochain $f$ and its image $\Psi
  (f)$\footnote{Only the non-zero values of $f$ appear on the figure.}}
\bigskip

\section{Three views of a secret}

 We give three equivalent constructions of the Steinberg
 representation. The fact that they are indeed equivalent is a
 consequence of the Borel-Serre theorem. 

 As usual, $H$ is the group of $F$-rational points of a split
 connected reductive $F$-algebraic group $\mathbb H$.

\subsection{The Steinberg representation via Zelevinski involution}

 We fix a maximal $F$-split torus
 $\mathbb T$ of $H$ as well as a Borel subgroup $\mathbb B$ containing
 $\mathbb T$. We denote by $\Phi$ the root system of $\mathbb H$
 relative to $\mathbb T$ (that we assume irreducible for simplicity
 sake),  by $\Phi^+$ the subset of positive roots
 relative to $\mathbb B$, and by $\Delta\subset \Phi^+$ the subset of
 simple roots. Recall\footnote{For more details, cf. \cite{Ca} {\S}1} that
 we have a bijection $\Theta \mapsto
 P_\Theta$ between the powerset of $\Delta$ and the set of parabolic
 subgroups of $G$ containing $B={\mathbb B}(F)$ (normalized by
 $P_\emptyset =B$ and $P_\Delta =H$). Each parabolic $P_\Theta$ has a
 standard Levi decomposition $P_\Theta =M_\Theta U_\Theta$, where
 $U_\Theta$ is (the group of $F$ rational points of) the unipotent
 radical of $P_\Theta$ and $M_\Theta$ a Levi component.  

 We denote by ${\mathcal R}(H)$ the category of smooth complex
 representations of $H$. For $\Theta\subset \Delta$,  ${\mathcal
   R}(M_\Theta )$ denotes the similar category attached to $M_\Theta$, 
 and ${\rm Ind}_{M_\Theta}^H$~: ${\mathcal R}(M_\Theta )\longrightarrow {\mathcal
   R}(H)$ , ${\rm Res}_H^{P_\Theta}$~: ${\mathcal R}(H)\longrightarrow {\mathcal
   R}(M_\Theta )$ the functors of normalized parabolic induction and normalized
 restriction (or Jacquet functor) respectively\footnote{Cf. \cite{Ca} {\S}3
   for more details.}. Both functors take representations of finite
 length to representations of finite length. In particular if $\rho$
 is a an irreducible smooth representations of $H$, the representation
 ${\rm Ind}_{P_\Theta}^G \, {\rm Res}_G^{P_\Theta}\, \rho$ gives rise
 to a well define element of $\displaystyle \left[{\rm Ind}_{P_\Theta}^G \, {\rm
     Res}_G^{P_\Theta}\, \rho\right]$ of the Grothendieck group
 ${\mathcal R}_{\rm fl}(H)$ of finite length smooth representations of
 $H$. 

 A.M. Aubert \cite{Au} has generalized the Zelevinski involution, defined
 by Zelevinsky for ${\rm GL}(N)$, to the case of any reductive
 group\footnote{Also see \cite{SS} for another construction of this
   involution.}.
 For an irreducible smooth representation $\rho$ of $H$, the
 Aubert-Zelevinski dual of $\rho$ is the element $\iota (\rho )$ of
 ${\mathcal R}_{\rm fl}(H)$ defined by
$$
\iota (\rho )=\sum_{\Theta\subset \Delta} (-1)^{\vert \Theta\vert} {\rm
    Ind}_{P_\Theta}^G \, {\rm Res}_G^{P_\Theta}\, \rho
$$
\noindent Then the key result of [Au] is that $\iota (\rho )$ is, up to a sign, an irreducible
representation of $H$. We shall denote this latter representation by
$\rho^\iota$. 

\smallskip

 \begin{definition} \label{Steinberg-Zelevinki} (The Steinberg representation via Zelevinski
   involution). One defines the Steinberg representation of $H$ to be
   the representation $({\mathbf 1}_H)^\iota$, that is the
   representation obtained by applying the Aubert-Zelevinski
   involution to the trivial representation  ${\mathbf 1}_H$ of $H$. 
\end{definition}

 It is a corollary of the proof of Borel-Serre theorem that the
 Steinberg representation of $H$, as previously defined, is in fact
 isomorphic to the top cohomology with compact support
 of the affine building of $H$ as an
 $H$-module. More precisely $({\mathbf 1}_H)^\iota$ is naturally
 isomorphic to the top reduced cohomology of the topological Tits
 building of $H$. 
\medskip

 In the case $H={\rm GL}(N,F)$, there is a simpler description of the
 Steinberg representation in terms of  parabolic induction. Take for
 $T={\mathbb T} (F)$ the subgroup of diagonal matrices and for $B$
 the Borel subgroup of upper triangular matrices. For $a\in F$, denote
 by $\vert a\vert_F$ absolute value of $a$ normalized by $\vert
 \varpi_F \vert_F = \frac{1}{q_F}$ for any uniformizer $\varpi_F$ of
 $F$. Finally let $\tau$ be the character of $T\simeq (F^\times )^N$
 given by
$$
\tau (t_1 ,...,t_N )= \vert t_1\vert_F^{(1-N)/2} \vert
t_2\vert_F^{(3-N)/2} \otimes\cdots \otimes  \vert t_N
\vert_F^{(N-1)/2}\ .
$$

\noindent Then the parabolically induced representation ${\rm Ind}_B^H \,
\tau$ has a unique irreducible $H$-quotient,  which turns out to be
the Steinberg representation of $H$. 

\medskip

Historically several definition of the Steinberg representation (or
of {\it special representations}) were given. That we give in this
section is the definition that Harish-Chandra gave in \cite{HC} {\S}15. In  \cite{Ca2}
Casselman proved that this representation, as defined by Harish-Chandra, is in fact irreducible.

\subsection{The Steinberg representation as a space of harmonic cochains}

 By the Borel-Serre theorem, the Steinberg representation ${{\mathbf S}{\mathbf t}}_H$ of
 $H$ is isomorphic to the top cohomology space $H_c^d (X_H ,{\mathbb C} )$ as
 a $H$-module. Since $X_H$ is a simplicial complex, it is a standard
 result of algebraic topology that the spaces $H_c^k (X_H ,{\mathbb C} )$ can
 be computed by {\it simplicial methods}. We recall the definition of
 the cohomological complex of alterned cochains on $X_H$ whose
 cohomology computes $H^*_c (X_H ,{\mathbb C} )$.

 Let $q\in \{ 0,...,d\}$. An {\it ordered $q$-simplex} of $X_H$ is an ordered
 sequence $(s_0 ,...,s_q )$ of vertices of $X_H$ such that $\{ s_0
 ,...,s_q\}$ is a $q$-simplex. We denote by $X_H^{(q)}$ the set of
 ordered $q$-simplices in $X_H$. The space $C_c^q (X_H )$ of
 alterned $q$-cochains on $X_H$ with finite support is the space of complex valued functions $f$~:
 $X_H^{(q)}\longrightarrow {\mathbb C}$ satisfying:
\smallskip

(a) $f$ has finite support,
\smallskip

(b) $f(s_{\tau (0)}, ..., s_{\tau (q)}) = {\rm sgn}(\tau )\, f(s_0
,...,s_q )$, for all ordered $q$-simplices $(s_0 ,...,s_q )$, all
permutations $\tau$ of the set $\{ 0,1,...,q\}$, where ${\rm sgn}$
denotes the signature of a permutation. 
\smallskip

Each $C_c^q (X_H )$ is endowed with a structure of smooth $H$-module
via the formula:
$$
(h. f)(s_0 ,...,s_q ) = f(h^{-1}.s_0 , ...,h^{-1}.s_q ), \ f\in C_c^q
(X_H ), \ h\in H, \ (s_0 ,...,s_q )\in X_H^{(q)}  \ .
$$
\noindent For $q=0, ..., d-1$, we defined a coboundary
map\footnote{The reader will forgive me to use the same symbol $d$ to
  denote the dimension of $X_H$ and the coboundary map.} $d$~: $C_c^q
(X_H )\longrightarrow C_c^{q+1}( X_H )$ by 
$$
(df)(s_0 ,...,s_{q+1}) = \sum_{i=0,...,q+1} (-1)^i f(s_0 , ...,
\hat{s}_i ,..., s_{q+1}), \  f\in C_c^q
(X_H ), \ (s_0 ,...,s_{q+1} )\in X_H^{(q+1)}\ ,
$$
\noindent where $(s_0 , ..., \hat{s}_i , ...,s_{q+1})$ denotes the ordered
$q$-simplex $(s_0 ,...,s_{i-1}, s_{i+1}, ..., s_{q+1})$. 
\medskip

We have a cohomological complex of smooth $H$-modules:
$$
0\longrightarrow C_c^0 (X_H )\overset{d}{\longrightarrow} \cdots \overset{d}{\longrightarrow} C_c^q
(X_H )\overset{d}{\longrightarrow} C_c^{q+1}(X_H )\overset{d}{\longrightarrow} \cdots
\overset{d}{\longrightarrow} C^d_c (X_H )\longrightarrow 0\ .
$$
\noindent Finally the $H$-modules $H_c^q (X_H ,{\mathbb C} )$, $q=0,...,d$,  are given by 
$$
H_c^q (X_H ,{\mathbb C} ) = {\rm ker} \left( d\ : \  C_c^q
(X_H ) \longrightarrow C_c^{q+1}(X_H )\right) / dC_c^{q-1}(X_H )
$$
\noindent with the convention that $C_c^i (X_H )=0$, for $i=-1,d+1$. In
particular we have ${{\mathbf S}{\mathbf t}}_H \simeq H_c^d (X_H ,{\mathbb C} )=C_c^d (X_H
)/dC_c^{d-1}(X_H )$. 
\medskip

We  now use the fact that $X_E$ is labellable in order to give a simpler
model of ${{\mathbf S}{\mathbf t}}_H$. So fix a labelling $\lambda$ of $X_E$ and for
$d$-simplex $C = \{ s_0 ,...,s_{d}\}$ of $X_H$, let
$C_\lambda$ be the unique ordered simplex $(t_0 ,...,t_d )$ such
that $\{ t_0 ,...,t_d \} =  \{ s_0 ,...,s_{d}\}$ and $\lambda (t_i
)=i$, $i=0,...,d$. Let ${\rm Ch}_H$ be the set of chambers of $X_H$ and
${\mathbb C}_c [{\rm Ch}_H ]$ be the set of complex valued functions on ${\rm Ch}_H$ with
finite support. Then we have an isomorphism of ${\mathbb C}$-vector spaces :
$$
{\mathbb C}_c [{\rm Ch}_H ]\longrightarrow C_c^d (X_H )
$$
\noindent given as follows. To any function $f\in {\mathbb C}_c [{\rm Ch}_H ]$, the
isomorphism attaches the unique $d$-cochain ${\tilde f}$ satisfying 
$$
{\tilde f}(C_\lambda ) = f(C)
$$
\noindent for any chamber $C$ of $X_H$. Because in general $H$ does not
preserves the orientation of chambers our isomorphism is not an
isomorphism of $H$-module in general (this is true if $\mathbb H$ is
simply connected). Define a $H$-module structure on ${\mathbb C}_c [{\rm Ch}_H ]$ by 
$$
(h. f)(C) =\epsilon_H (h)\, f(h^{-1}C), \ C \ {\rm chamber} \ {\rm
  of}\ X_H\ ,
$$
\noindent where $\epsilon_H$ is the quadratic character of $H$ defined in
{\S}1.3. Write ${\mathbb C}_c [{\rm Ch}_H ]\otimes \epsilon_H$ for the ${\mathbb C}$-vector
space ${\mathbb C}_c [{\rm Ch}_H ]$ endowed with the $H$-module structure we have just
defined. We then have :
\begin{lemma}
The map $f\mapsto {\tilde f}$ induced an isomorphism of smooth
$H$-modules: $ {\mathbb C}_c [{\rm Ch}_H ]\otimes
\epsilon_H \longrightarrow C_c^d (X_H )$. 
\end{lemma}

The $H$-module $dC_c^{d-1} (X_H )$ is generated by the functions
$df_D$, where $D= \{ t_0 ,...,t_{d-1}\}$ runs over
the $(d-1)$-simplices of $X_H$, and where $f_{(t_0 ,...,t_{d-1})}$ is
a $(d-1)$-cochain with support the set of ordered $(d-1)$-simplices of
the form $(t_{\tau (0 )}, ..., t_{\tau (d-1)})$, where $\tau$ is any
permutation of $\{ 0,...,d-1\}$. It is an easy exercise to show that
for all codimension $1$ simplex $D$,
up to a sign, and through the isomorphism $C_c^d (X_H )\simeq {\mathbb
  C}_c [{\rm Ch}_H ]$,  $df_D$ is the characteristic function of the set of chambers $C$
containing $D$. So we have proved:

\begin{proposition} As an $H$-module, ${{\mathbf S}{\mathbf t}}_H$ is the quotient of ${\mathbb C}_c
  [{\rm Ch}_H ]\otimes \epsilon_H$ by the submodule ${\mathbb C}_c
  [{\rm Ch}_H ]^0$ spanned by
  the function of the form $g_D$, $D$ codimension $1$ simplex of
  $X_H$, where for all chambers $C$
$$
g_D (C)=
\left\{ \begin{array}{ccl}
1 & {\rm if} & C\supset D\\
0 & {\rm ortherwise} &
\end{array}\right.
$$
\end{proposition}

 Recall that, being given a smooth representation $(\pi ,{\mathcal V})$ of $H$,
 we have two notions of {\it dual representations}. The {\it algebraic
   dual} is the representation $(\pi^* , {\mathcal V}^* )$ where ${\mathcal V}^* = {\rm
   Hom}_{\mathbb C} \, ({\mathcal V}, {\mathbb C} )$ is the
 space of linear forms on ${\mathcal V}$, and where $H$ acts by 
$$
(\pi^* (h)\varphi )(v) = \varphi (\pi (h^{-1})v), \ h\in H, \ \varphi
 \in {\mathcal V}^* , \ v\in {\mathcal V} \ .
 $$
 
\noindent The {\it smooth dual} or {\it contragredient} is the
sub-$H$-module $({\tilde \pi}, {\tilde V})$ of  $(\pi^* , {\mathcal V}^* )$
formed of {\it smooth} linear forms, that is linear forms fixed by an
open  subgroup of $H$. We are going to give very simple models
for ${{\mathbf S}{\mathbf t}}^*$ and
${\tilde {{\mathbf S}{\mathbf t}}}_H\subset {{\mathbf S}{\mathbf
    t}}_H^*$.
We shall see in the next section that the  representation
${{\mathbf S}{\mathbf t}}_H$ is self-dual. So we shall obtain a simple model for ${{\mathbf S}{\mathbf t}}_H$ as
well. 

In this aim, observe that there is a perfect pairing $\langle
-,-\rangle$ bewteen the $H$-modules ${\mathbb C}_c [{\rm Ch}_H ]$
and ${\mathbb C} [{\rm Ch}_H ]$, the
space of complex valued functions on $X_H$ with no condition of
support. It is given by 
$$
\langle f, \omega \rangle =\sum_{C\in {\rm Ch}_H} f(C)\omega (C), \ f\in
{\mathbb C}_c [X_H ], \ \omega \in {\mathbb C} [X_H ]\ .
$$
\noindent Define the space of {\it harmonic $d$-cochains} on $X_H$ as the
orthogonal ${\mathcal H} (X_H )$ of ${\mathbb C}_c [{\rm Ch}_H ]^0$ in
${\mathbb C} [{\rm Ch}_H]$ relative to
the pairing $\langle -,-\rangle$. In other words, an element
$\omega\in {\mathbb C} [{\rm Ch}_H ]$ is said to be {\it harmonic} if it satisfies
the {\it harmonicity condition}:
$$
\sum_{C\supset D} \omega (C)=0
$$
\noindent for all codimension $1$ simplex $D$ of $X_H$. Finally define
${\mathcal H} (X_H )^\infty$ to be the space of smooth  harmonic
$d$-cochains, that is harmonic $d$-cochains which are fixed by an open
subgroup of $H$. As a consequence of the previous proposition we have
the following. 

\begin{proposition}
We have two isomorphisms of $H$-modules :
$$
{{\mathbf S}{\mathbf t}}_H^* =
{\rm Hom}_{\mathbb C} ({{\mathbf S}{\mathbf t}}_H ,{\mathbb C} )\simeq {\mathcal H} (X_H )\otimes
\epsilon_H , \ \ {{\mathbf S}{\mathbf t}}_H \simeq
\tilde {{\mathbf S}{\mathbf t}}_H \simeq {\mathcal H} (X_H )^\infty
\otimes \epsilon_H\ .
$$
\end{proposition}

\medskip

Smooth harmonic cochains on $X_H$  are quite tricky objects. One can for instance
prove that the unique smooth harmonic cochain with finite support
 is the zero cochain. 
 In the next next section, we shall exhibit a non-trivial smooth
 harmonic cochain: an Iwahori spherical vector in ${\mathcal H} (X_H )^\infty$.

\subsection{The Steinberg representation via Type Theory}

\medskip

In this section we assume that the algebraic group $\mathbb H$ is
simply connected (e.g. ${\mathbb S}{\mathbb L}_n$, ${\mathbb S}{\rm
  p}_{2n}$). Fix an Iwahori subgroup $I$ of $H$ and consider the full
subcategory ${\mathcal R}_I (H)$ of ${\mathcal R}(H)$  defined as follows:
a smooth representation $(\pi ,{\mathcal V} )$ of $H$ in a ${\mathbb C}$-vector space
${\mathcal V}$ is an object of ${\mathcal R}_I (H)$ if, as a $H$-module, ${\mathcal V}$
is generated by the subset ${\mathcal V}^I$ of vectors fixed by $I$. 
\smallskip

 By \cite{Bo}  this category may be described in terms of parabolic
 induction. One says that an irreducible smooth representation $\pi$ of $H$
 belongs to the {\it unramified principal series} if there exists an
 unramified character $\chi$ of a maximal split torus $T$ of $H$ such
 that $\pi$ is a subquotient of ${\rm Ind}_B^H\, \chi$, for some Borel
 subgroup $B$ containing $T$. Here {\it unramified} means that $\chi$
 is trivial on the maximal compact subgroup of $T$. Then a
 representation $(\pi ,{\mathcal V} )$ is an object of ${\mathcal R}_I (H)$ if
 and only if all irreducible subquotients of $\pi$ belong to the
 unramified principal series. 
\smallskip

 In particular the category ${\mathcal R}_I (H)$ is stable by the
 operation of taking subquotient. In the terminology of Bushnell and
 Kutzo's {\it theory of  types} (cf. \cite{BK} for a fondation of this
 theory),  one says that the pair $(I, {\mathbf
   1}_I )$ is a {\it type} for $H$. 
\smallskip

 Let $\mu$ be a Haar measure on $H$ normalized by $\mu (I)=1$. Let
 ${\mathcal H}(H)$ be the space of complex locally constant functions
 on $H$ with compact support. Let
 ${\mathcal H}(H,I)$ be the ${\mathbb C}$-vector space of bi-$I$-invariant
 complex functions on $H$ with compact support. Equip ${\mathcal
   H}(H)$ and ${\mathcal H}(H,I)$ with
 the convolution product:
$$
f_1 \star f_2 (h)=\int_H f_1 (hx)f_2 (x^{-1})\, d\mu (x)\ .
$$
Then ${\mathcal H}(H)$ is an associative algebra and ${\mathcal H}(H,I)$ is a
subalgebra  with unit $e_I$, the
characteristic function of $I$, called the {\it Iwahori--Hecke}
algebra of $H$. In fact ${\mathcal H}(H,I) =e_I \star {\mathcal H}(H)\star
e_I$. This latter algebra  is non  commutative if the semisimple rank of $H$
is $>0$. Recall that is it a basic fact of the theory of smooth
representations of $p$-adic reductive groups that the categories
${\mathcal R}(H)$ and ${\mathcal H}(H)-{\rm Mod}$ (the category of
{\it non degenerate}\footnote{ An ${\mathcal H}(H)$-module $M$ is non
  degenerate if ${\mathcal H}(H)\cdot M =M$.} left ${\mathcal H}(H)$-modules) 
  are ``naturally'' isomorphic
(cf. \cite{Ca}  for more details). 

 If $(\pi , {\mathcal V} )$ is  smooth representation of $H$, then ${\mathcal V}^I$ is
 naturally a left  ${\mathcal H}(H,I)$-module. In particular we have a functor
 $m_I$~: ${\mathcal R}_I (H)\longrightarrow {\mathcal H}(H,I) - {\rm Mod}$, defined by 
$(\pi ,{\mathcal V} )\mapsto {\mathcal V}^I$. Historically the following result is the
 keystone of Type Theory.

\begin{theorem} \label{EqCat}(Cf. \cite{Bo}, \cite{BK})  The functor $m_I$ is an equivalence
  of categories. An inverse $M_I$ of $m_I$ is given by
$$
M_I \ : \ {\mathcal H}(H,I) - {\rm Mod}\longrightarrow {\mathcal R}_I (H)\ , M\mapsto
{\mathcal H}(H)\otimes_{{\mathcal H}(H,I)} M\ .
$$
\end{theorem}

\noindent In the theorem the $H$-module structure of ${\mathcal H}(H)\otimes_{{\mathcal H}(H,I)} M$
comes from the action of $H$ on ${\mathcal H}(H)$ by left
translation. 
\medskip

 Fix a maximal split torus $\mathbb T$ of $\mathbb H$ such that the
 chamber $C$ fixed by $I$ lies in the apartment ${\mathcal A}$ attached to $\mathbb
 T$. Recall that we have the Iwahori decomposition 
$$
H= IW^{\rm Aff} I=\bigsqcup_{w\in W^{\rm Aff}} IwI
$$
\noindent where $W^{\rm Aff}$ is the affine Weyl group attache to $\mathbb
T$. Also recall that the  Coxeter group $W^{\rm Aff}$ is generated by a finite set of
involutions $S$ (attached to the pair $(C, {\mathcal A} )$. 

It follows from the Bruhat-Iwahori decomposition, that as  a ${\mathbb C}$-vector
space,  ${\mathcal H}(H,I)$ has basis $(e_w)_{w\in W^{\rm Aff}}$, where $e_w$ is
the characteristic function of $IwI$. 
The structure of the algebra ${\mathcal H}(H,I)$ is well known.

\begin{theorem} (Iwahori-Matsumoto \cite{IM}) The unital ${\mathbb C}$-algebra
  ${\mathcal H}(H,I)$ has the following presentation: it is generated by the
  $e_s$, $s\in S$, with the relations
\smallskip

(R1) $e_s^2 = (q_F -1)e_s +q_F e_1$, $s\in S$,
\smallskip

(R2) for all distinct $s$, $t$ in $S$, we have
\smallskip

\quad $(e_s e_t)^r  e_s = e_t (e_s e_t )^r$, if $m_{st}=2r+1$,
\smallskip

\quad $(e_s e_t )^r =(e_t e_s )^r$, if $m_{st}=2r$,
\smallskip

\noindent where $m_{st}$ is the order of $st\in W^{\rm Aff}$. 
\end{theorem}

The quadratic relations (R1) writes $(e_s +1)(e_s -q_F )=0$, $s\in S$.  It follows
that the algebra ${\mathcal H}(H,I)$ admits a unique character $\chi$
(equivalently a $1$-dimensional left module) defined by $\chi (e_s
)=-1$.  This character is known as the {\it special character} of
${\mathcal H}(H,I)$. By the equivalence of categories \ref{EqCat}, $\chi$
corresponds to an irreducible  smooth representation $(\pi_\chi ,{\mathcal V}_\chi
)$ of $H$. We are going to prove that this representation is nothing other
than the Steinberg representation of $H$.
\medskip

We have ${\mathcal V}_\chi ={\mathcal H}(H)\otimes_{{\mathcal H}(H,I)} {\mathbb C}$ where ${\mathcal H}(H,I)$
acts on ${\mathbb C}$ via the character $\chi$. Since $e_I$ is the unit
element of ${\mathcal H}(H,I)$, this  may be rewritten $\displaystyle {\mathcal V}_\chi = {\mathcal
  H}(H)\star e_I \otimes_{{\mathcal H}(H,I)} {\mathbb C}$. The $(H,{\mathcal H}(H,I) )$-bimodule
${\mathcal H}(H)\star e_I$ is the space of locally constant function on
$G$ which are right $I$-invariant and have compact support. Since $I$
is the global stabilizer of a chamber of $X_H$, the discrete
topological space
$H/I$ is isomorphic to the set of chambers in $X_H$ as a $H$-set; we
denote by ${\rm Ch}_H$ this set of chambers.  It
follows that ${\mathcal H}(H)\star e_I$ identifies with ${\mathbb C}_c [{\rm Ch}_H ]$, the
set of complex valued functions with finite support on ${\rm Ch}_H $. Under
this identification, the left $H$-module structure of ${\mathbb C}_c [{\rm Ch}_H ]$ is
the natural one. 

 The Bruhat-Iwahori decomposition $I\backslash H/I \simeq W^{\rm Aff}$
 allows us to classify the relative positions of two chambers of
 $X_H$, that is the orbits of $H$ in ${\rm Ch}_H\times {\rm Ch}_H$ : two
 chambers $C_1$, $C_2$ are in position $w\in W^{\rm Aff}$, denoted by 
$C_1 \sim_w C_2$,  if the $H$-orbit of $(C_1 ,C_2 )$ contains $(C_0 ,
 wC_0 )$, where $C_0$ is the chamber fixed by $I$. The following lemma
 is an excellent exercise left to the reader. 

\begin{lemma} Under the natural identification ${\mathcal H}(H)\star
  e_I\simeq {\mathbb C}_c [{\rm Ch}_H ]$, the right ${\mathcal H}(H,I)$-module structure of 
${\mathbb C}_c [{\rm Ch}_H ]$ is given as follows:
$$
f\star e_s (C) = \sum_{C'\sim_s C} f(C') = \sum_{C'\supset C_s ,
  \ C'\not= C}f(C'), \ f\in {\mathbb C} [{\rm Ch}_H ], \ s\in S, \ C\in {\rm Ch}_H \ ,
$$
\noindent where $C_s$ denotes the codimension $1$ subsimplex of $C$ of type
$s$. 
\end{lemma}

By definition the tensor product 
$$
 {\mathcal  H}\star e_I \otimes_{{\mathcal H}(H,I)} {\mathbb C}
\simeq {\mathbb C}_c [{\rm Ch}_H ]\otimes_{{\mathcal H}(H,I)}{\mathbb C}
$$
\noindent is the quotient of ${\mathbb C}_c [{\rm Ch}_H ]\otimes_{\mathbb C} {\mathbb C}
= {\mathbb C}_c [{\rm Ch}_H ]$ by the subspace  generated by the
functions  $f\star e_w -\chi(w) f$,
where $f$ runs over ${\rm Ch}_H$ and $w$ runs over $W^{\rm Aff}$. Since
the $e_s$, $s\in s$, generate ${\mathcal H}(H,I)$ as an algebra, this
subspace is also generated by the $f\star e_s -\chi (s) f = f\star e_s
+f$. By the previous lemma,  this is the space of functions
generated by the  $f_D$, $D$ codimension $1$ simplex of $X_H$, defined
by $f_D (C)=1$ if $C\supset D$, $f_D (C)=0$ otherwise. So this space
is nothing other than the space ${\mathbb C}_c [{\rm Ch}_H ]^0$ defined in {\S}3.2. 
The following is now a consequence of Proposition 3.3.

\begin{proposition} The $H$-module  ${\mathcal V}_\chi ={\mathcal
    H}(H)\otimes_{{\mathcal H}(H,I)} {\mathbb C}$
is isomorphic to ${{\mathbf S}{\mathbf t}}_H \simeq H_c^d (X_H ,{\mathbb C} )$. 
\end{proposition}

\begin{corollary} a) The Steinberg representation  of $H$ has non-zero
  fixed vectors under the Iwahori subgroup $I$. Moreover ${{\mathbf S}{\mathbf t}}_H^I$ is
  $1$-dimensional.

\noindent b) The Steinberg representation is self-dual. 
\end{corollary}

Only b) needs to be proved. The ${\mathcal H} (H,I)$-module $m_I ({\tilde
  {{\mathbf S}{\mathbf t}}}_H) = ({\tilde {{\mathbf S}{\mathbf t}}}_H )^I$ is the dual of ${{\mathbf S}{\mathbf t}}_H^I$, and this
latter module is $1$-dimension, the algebra ${\mathcal H} (H,I)$ acting via the
character $\chi$. Since $\chi$ has real values, the ${\mathcal H} (H,I)$-module
${{\mathbf S}{\mathbf t}}_H^I$ is self-dual. It follows that ${{\mathbf S}{\mathbf t}}_H$ is self dual since
$m_I$ is an equivalence of categories. 
\medskip

A non-zero vector vector in ${{\mathbf S}{\mathbf t}}_H^I$ is called {\it
  Iwahori-spherical}. In the next proposition, we describe the line
${{\mathbf S}{\mathbf t}}_H^I$ in the model ${{\mathbf S}{\mathbf
    t}}_H \simeq  {\mathcal H} (X_H )^\infty$. This will
exhibit a non-trivial element of ${\mathcal H} (X_H )^\infty$. 

\begin{proposition} Let $C_0$ denote the chamber fixed by $I$. 
The exists a unique Iwahori-spherical vector $f_{C_0}$
  in   $ {\mathcal H} (X_H )^\infty$ satisfying $f_{C_0}(C_0 )=1$. It is given
  by
$$
f_{C_0}(C) = \left(\frac{-1}{q_F}\right)^{d(C_0 ,C)}, \ C\in {\rm Ch}_H\  .
$$
\noindent In particular if $C$ is a chamber of the apartment attached to
the torus $\mathbb T$, we have
$$
f_{C_0}(C)=  \left(\frac{-1}{q_F}\right)^{l(w)}, \ {\rm if}\  C=wC_0 ,
\ w\in W^{\rm Aff}\  .
$$
\end{proposition}

 Indeed let us first remark that $f_{C_0}$ is $I$-invariant; this is
 due to the fact that, since the action of $H$ on $X_H$ is simplicial,
 it preserves the distance $d$ between pairs of chambers. In particular,
 $I$ being open, $f_{C_0}$ is a smooth function. Let us prove that it
 is harmonic. Let $D$ be a codimension $1$ chamber. We need the
 following lemma whose proof we shall admit. 

\begin{lemma} There exists a unique chamber $C_1$ in $X_H$ containing
  $D$ and  such that the
  distance $\delta =d(C_0 ,C_1)$ is minimal. In particular there
  exists an integer $\delta \geqslant 0$ such that among the
  $q_F +1$ chambers containing $D$, one is at distance $\delta$ from
  $C_0$ and the other at distance $\delta +1$. 
\end{lemma}

\medskip

Let $C_1$ and $\delta$ be as in the lemma. We have:
\begin{align*}
\sum_{C\supset D} f_{C_0}(C) &= \sum_{C\supset D} \left(\frac{-1}{q_F}\right)^{d(C_0 ,C)} \\
 & =
\left(\frac{-1}{q_F}\right)^\delta
+ q_F \left(\frac{-1}{q_F}\right)^{\delta +1} \\
 &= 0
\end{align*}
\noindent so that the harmonicity condition at $D$ holds true. 
 
\section{Distinction of the Steinberg representation} We fix a Galois
symmetric space $G_E /G_F$ as in the introduction. So $E/F$ is a Galois
quadratic extension of non-archimedean local fields and we have $G_E
={\mathbb G}(E)$, $G_F ={\mathbb G}(F)$, where $\mathbb G$ is a
connected reductive group defined over $F$. In \cite{Pr}{\S}7,  assuming
that  the derived group ${\mathbb G}^{\rm der}$ is quasi-split over
$F$, D. Prasad defines a quadratic
character $\epsilon_{\rm Prasad}$ of $G_F$ and makes the following
conjecture.

\begin{conjecture} Assume that ${\mathbb G}^{\rm der}$ is quasi-split and let
  ${{\mathbf S}{\mathbf t}}_E$ denote the Steinberg representation of $G_E$.
  \smallskip

  (1) The intertwing space ${\rm Hom}_{G_F} ({{\mathbf S}{\mathbf t}}_E ,\epsilon_{\rm
    Prasad})$ is $1$-dimensional.
  \smallskip

  (2) For any character $\chi$ of $G_F$ such that $\chi\not=
  \epsilon_{\rm Prasad}$, we have  ${\rm Hom}_{G_F} ({{\mathbf S}{\mathbf t}}_E ,\epsilon_{\rm
    Prasad}) = 0$.
\end{conjecture}

As explained in the introduction,  this statement is in fact a particular case of a much more general
 conjecture of Prasad's which predicts the distinction of an
 irreducible representation of $G_E$ in terms of its Galois parameter 
 (that is through the conjectural local  Langlands correspondence).
 \medskip

 Conjecture 4.1 is proved in \cite{BC} under the following assumptions:
 \medskip

 (H1) $\mathbb G$ is
 split over $F$,

 (H2) the  adjoint group of $\mathbb G$ is simple,

 (H3) the extension
 $E/F$ is unramified.
 \medskip
 
 \noindent In fact Assumption (H2) can easily be removed as shown in \cite{Cou}.  In this
  section we assume that (H1), (H2), (H3) hold. We give some hints for
  the proof provided in \cite{BC} and
  we make it simpler by the use of Poincar\'e series.  In {\S}4.5 we
  shall say a few words on this extension of this result, extension due to
  Fran\c cois Court\`es, to the case where $E/F$ is {\it tamely
    ramified}. 

\subsection{The invariant linear form} As in {\S}1.3, we denote by
$X_F$ (resp. $X_E$) the  semisimple Bruhat-Tits building of $G_F$
(resp. of $G_E$). Since $E/F$ is unramified, we have a natural
embedding $X_F \subset X_E$ which is simplicial, $G_F$-equivariant and
${\rm Gal}(E/F)$ equivariant. In particular the set ${\rm Ch}_F$ of
chambers of $X_F$ is naturally a subset of ${\rm Ch}_E$, the set of
chambers of $X_E$.

By Proposition 3.4, the Steinberg representation of $G_E$ is given by
${{\mathbf S}{\mathbf t}}_E \simeq {\mathcal H}(X_E )^\infty \otimes \epsilon_E$, where:
\smallskip

 -- ${\mathcal H}(X_E )$ is,  as defined in {\S}3.2, the space of harmonic $d$-cochains on
 $X_E$ ($d$ is here the semisimple rank of $\mathbb G$), and
 ${\mathcal H}(X_E )^\infty$ is the subspace of $G_E$-smooth vectors.

  -- $\epsilon_E = \epsilon_{G_E}$ is the quadratic character of $G_E$
  defined in {\S}1.2.
  \smallskip

  It turns out \cite{Cou} that the restriction $\epsilon_{\vert G_F}$
  coincides with Prasad's character $\epsilon_{\rm Prasad}$. It
  follows that the intertwining space ${\rm Hom}_{G_F} ({{\mathbf S}{\mathbf t}}_E
  ,\epsilon_{\rm Prasad})$ is given by
  $$
  {\rm Hom}_{G_F}({\mathcal H}(X_E )^\infty \otimes \epsilon_E , \epsilon_{\rm
    Prasad}) = {\rm Hom}_{G_F}({\mathcal H}(X_E)^\infty , {\bf 1}),
  $$

  \noindent where ${\bf 1}$ denotes the trivial character of $G_F$.
  \medskip

   So in order to prove Conjecture 4.1 in our case, we have to establish: 
   \smallskip

   (1) $\displaystyle {\rm dim}\,  {\rm Hom}_{G_F^{\rm der}}({\mathcal H}(X_E)^\infty ,
   {\bf 1})\leqslant 1$,

    (2) $ {\rm Hom}_{G_F}({\mathcal H}(X_E)^\infty , {\bf 1}) \not=
   0$,
   \smallskip

   \noindent where $G_F^{\rm der}$ denotes the derived group of $G$.
   \smallskip
   
    The proofs of (1) and (2) are quite different in nature. We shall 
    say a few words on   the proof of (1) in {\S}4.4 and  we refer to
    [BC] for more details. To
    prove (2) we have to exhibit a non-zero $G_F$-invariant linear
    form
    $$
    \Lambda~: \ {\mathcal H}(X_E)^\infty \longrightarrow {\mathbb C} \ .
    $$

    \noindent It is quite natural to set
    \begin{equation}
      \Lambda (f) =\sum_{C\in {\rm Ch}_F} f(C) , \ f\in {\mathcal H}(X_E
      )^\infty
      \end{equation}
\noindent since, if $\Lambda$ is well defined, it is clearly linear and
$G_F$-equivariant. 
\smallskip

    \noindent Of course we have to prove that for each $f\in {\mathcal H}(X_E
      )^\infty$ the sum of (4) converges and that there exists $f_0
    \in {\mathcal H}(X_E )^\infty$ such that $\Lambda (f_0 )\not= 0$
    (such a vector $f_0$ is called a {\it test vector} for
    $\Lambda$). More precisely we prove the following.

    \begin{proposition} (1) If $f\in {\mathcal H}(X_E )^\infty$, then
      the restriction $f_{{\rm Ch}_F }$ lies in $L^1 ({\rm Ch}_F )$, the space
      of summable complex functions on ${\rm Ch}_F$.

      (2) We have $\Lambda (f_{\rm Iwahori}) \not= 0$, where $f_{\rm
        Iwahori}$ is the Iwahori-spherical vector relative to some
      fixed chamber $C$ of $X_F$.
    \end{proposition}

    In (4.3) we shall give a proof of this proposition which differs
    from that of \cite{BC} (and which is much simpler). It relies on a good
    understanding of the combinatorics of chambers in $X_F$ thanks to
    the use of Poincar\'e series. 
    
\subsection{Combinatorics of chambers}  We fix a maximal $F$-split torus
$\mathbb T$ of $\mathbb G$. Let $T={\mathbb T}(F)$ and $N_G (T)$ be
the normalizer of $T$ in $G$. The spherical Weyl group $N_G (T)/T$ is
denoted $W^{\rm Sph}$.    Let ${\mathcal A}$ be the
apartment of $X_F$ attached to $T$; this is also the (Galois fixed)
apartment of $X_E$ attached to  ${\mathbb T}(E)$. Fix a chamber $C$ in
${\mathcal A}$ and write $I$ for the Iwahori subgroup of $G_F$ fixing $C$. Let
$W^{\rm Aff} = N(T)/(T\cap I)$ be the extended affine Weyl group of
$G_F$. As in {\S}1.1, this Weyl group decomposes as $W^{\rm Aff} =
\Omega\rtimes W^{\rm Aff}_0$, where $W^{\rm Aff}_0$ is an affine
Coxeter group. We denote by $l$ the length function on $W^{\rm Aff}_0$
attached to the chamber $C$.
\medskip

  The group $G_0 := IW_0^{\rm Aff}I$ is a (normal) subgroup of $G$ which
  acts transitively on  ${\rm Ch}_F$. So we may write the disjoint union
  decomposition:
  \begin{equation}
  {\rm Ch}_F = \bigsqcup_{w\in W_0^{\rm Aff}} \bigsqcup_{g\in IwI/I} \{ g.C\},
  \end{equation}

\noindent where the fact that the union is indeed disjoint comes from the
  fact that ${\mathcal A}\cap {\rm Ch}_F$ is a fundamental domain for the action of
  $I$ on ${\rm Ch}_F$. For future calculations, we need a formula for the
  cardinal of $IwI/I$, $w\in W_0^{\rm Aff}$.

  \begin{proposition} (\cite{IM} Prop. 3.2) For $w\in W_0^{\rm Aff}$, we
    have:
    $$
    \left\vert IwI/I \right\vert = q_F^{l(w)}\ ,
    $$

  \noindent where $q_F$ is the cardinal of the residue field of $F$. 
\end{proposition}

  Let $N(d)$ denote the number of chambers in ${\mathcal A}$ at combinatorial
  distance $d$ from $C$. By definition the Poincar\'e series of
  $W_0^{\rm Aff}$ is the generating function
  $$
  P_{W_0^{\rm Aff}} (X) = \sum_{k\geqslant 0} N(k) X^k = \sum_{w\in W_0^{\rm Aff}}
  X^{l(w)}\ .
  $$

  A close formula for  this Poincar\'e series is known:

  \begin{theorem} (\cite{Bott}, \cite{St}) The formal series $P_{W_0^{\rm Aff}}$
    is a rational function given by
    $$
    P_{W_0^{\rm Aff}} (X) = \frac{1}{(1-X)^{d-1}}\prod_{i=1}^{d-1}
      \frac{1-X^{m_i}}{1-X^{m_i -1}}
    $$
    \noindent where $m_1$, $m_2$, ..., $m_{d_1}$ are the exponents of the
    finite Coxeter group $W^{\rm Sph}$ (see \cite{Bou} Chap. V, {\S}6,
    D\'efinition 2).

    In particular, the radius of convergence of  $P_{W_0^{\rm Aff}} (X)$
is $1$ and $P_{W_0^{\rm Aff}}$ defines a non-vanishing function on the
real open  interval $(-1,1)$.
\end{theorem}

  For instance if $W^{\rm Sph}$ is of type $A_l$ (case of ${\rm
    GL}_{l+1}$ or ${\rm SL}_{l+1}$), then we have $m_i = i$,
  $i=1,...,l$ (cf. \cite{Bou} Planche I). 

\subsection{The Poincar\'e series trick} We begin by  proving that for $f\in
{\mathcal H} (X_E)^\infty$, the infinite sum (4) defining $\Lambda (f)$ is absolutely convergent,
 that is $f_{\vert {\rm Ch}_F}\in L^1 ({\rm Ch}_F )$. For this we first use the fact that if a function
$f$~: ${\rm Ch}_E\longrightarrow {\mathbb C}$ statisfies the hamonicity condition and is smooth under
 the action of $G_E$, then it decreases in a way described as follows.

\begin{proposition} (Cf. \cite{BC})  Let $f\in {\mathcal H} (X_E )^\infty$. There exists a real 
$K_f >0$ such
 that for all chambre $D$ of $X_E$, we have
$$
\vert f(D) \vert \leqslant K_f\, q_E^{-d(C,D)} \ ,
$$
\noindent where $q_E = q_F^2$ is the cardinal of the residue field $k_E$, and where $d(C,D)$
denotes the combinatorial distance between chambers of $X_E$. 
\end{proposition}

Now, for $f\in {\mathcal H} (X_E )^\infty$, using decomposition (5), we may write:
\begin{align*}
\sum_{D\in {\rm Ch}_F} \vert f(D)\vert & \leqslant K_f \sum_{D\in {\rm Ch}_F} q_E^{-d(C,D)}\\
 & \leqslant K_f \sum_{w\in W_0^{\rm Aff}} \sum_{g\in IwI/I} q_E^{-d(C,gC)}
\end{align*}

\noindent If $g\in IwI$ for some $w\in W_0^{\rm Aff}$, we write
 $g=i_1 w i_2$, with $i_1$, $i_2\in I$, so that 
$$
d(C,gC)=d(i_1^{-1}C,wi_2 C)=d(C,wC)=l(w)
$$
\noindent where we used the facts that the distance $d$ is $G_E$-invariant and that 
$C$ is fixed by $I$. So we obtain:

\begin{align*}
\sum_{D\in {\rm Ch}_F} \vert f(D)\vert & \leqslant K_f \sum_{w\in W_0^{\rm Aff}} 
\sum_{g\in IwI/I} q_E^{-l(w)}\\
& \leqslant K_f \sum_{w\in W_0^{\rm Aff}} \vert IwI/I\vert q_E^{-l(w)}\\
& \leqslant K_f  \sum_{w\in W_0^{\rm Aff}} q_F^{l(w)} q_E^{-l(w)}\\
& \leqslant K_f  \sum_{w\in W_0^{\rm Aff}} (\frac{1}{q_F})^{l(w)}\\
& \leqslant K_f P_{W_0^{\rm Aff}} (\frac{1}{q_F})
\end{align*}

\noindent where we used the fact that $\vert IwI/I\vert =q_F^{l(w)}$ (Proposition 4.3) and that 
$q_E = q_F^2$. Now since the radius of convergence of the series $P_{W_0^{\rm Aff}}$ is $1$, we
obtain  $\displaystyle P_{W_0^{\rm Aff}}(\frac{1}{q_F})<+\infty$ and the sum defining $\Lambda (f)$
is indeed convergent. 
\medskip

We now prove that $\Lambda$ is non-zero by computing its value at the Iwahori fixed vector
of ${{\mathbf S}{\mathbf t}}_E$ given in {\S}3.3. Recall that it is given by
$$
f_{\rm Iwahori}(D)=\left(\frac{-1}{q_E}\right)^{d(C,D)}, \ \ D\in {\rm Ch}_E\ .
$$
We have:
\begin{align*}
\Lambda (f) & = \sum_{D\in {\rm Ch}_F} (\frac{-1}{q_E})^{d(C,D)}\\
& = \sum_{w\in W_0^{\rm Aff}} \sum_{g\in IwI/I} (\frac{-1}{q_E})^{d(C,D)}\\
& = \sum_{w\in W_0^{\rm Aff}} \vert IwI/I\vert (\frac{-1}{q_F^2})^{l(w)}\\
& = \sum_{w\in W_0^{\rm Aff}} q_F^{l(w)} (\frac{-1}{q_F^2})^{l(w)}\\
& = P_{\rm W_0^{\rm Aff}} (-\frac{1}{q_F})
\end{align*}

Since $P_{\rm W_0^{\rm Aff}}$ does not vanish on the open interval
$(-1,1)$,  we have proved the following result. 

\begin{proposition} Let $f_{\rm Iwahori}\in {{\mathbf S}{\mathbf t}}_E$
  be a non-zero Iwahori-spherical vector
  of ${{\mathbf S}{\mathbf t}}_E$ and $\Lambda \in
  {\rm Hom}_{G_F} ({{\mathbf S}{\mathbf t}}_E ,\epsilon_{\rm Prasad})$ be a non zero
equivariant linear form. Then $\Lambda (f_{\rm Iwahori})\not= 0$, and for suitable 
normalizations of $f_{\rm Iwahori}$ and $\Lambda$, we have the formula:
$$
\Lambda (f_{\rm Iwahori})=P_{W_0^{\rm Aff}}(-\frac{1}{q_F}) = 
\frac{1}{(1+\frac{1}{q_F})^{d-1}}\prod_{i=1}^{d-1}
      \frac{1-(-\frac{1}{q_F})^{m_i}}{1-(-\frac{1}{q_F})^{m_i -1}}
$$
\noindent where $d$ is the rank of the spherical Weyl group $W^{\rm Sph}$ of $\mathbb G$
and $m_1$, ...,$m_d$ the exponants of $W^{\rm Sph}$. 
\end{proposition}

Of course, once one knows that $\Lambda (f_{\rm Iwahori})\not= 0$,  one can always find
normalizations so that the previous formula holds. The point is that such normalizations
are natural in the model of ${{\mathbf S}{\mathbf t}}_E$ given by smooth harmonic cochains.

\subsection{ Multiplicity one} The proof of
the multiplicity $1$ property, i.e. assertion (1) of {\S}4.1,  proceeds
as follows.  We use the natural $G_E^{\rm der}$-isomorphism ${\rm Hom}_{\mathbb C}
({{\mathbf S}{\mathbf t}}_E ,{\bf 1}) \simeq {\mathcal H} (X_ E)$ so that (1) may be rewritten:
\smallskip

(1') \ \ ${\rm dim}\, {\mathcal H} (X_E )^{G_F^{\rm der}} \leqslant 1$,
\smallskip

\noindent where ${\mathcal H} (X_E )^{G_F^{\rm der}}$ denotes the ${\mathbb C}$-vector space
of $G^{\rm der}$-invariant harmonic cochains. Let us fix a chamber
$C_0$ in $X_F$. The basic idea is to
prove that the map
$$
j~: \ {\mathcal H} (X_E )^{G_F^{\rm der}} \longrightarrow {\mathbb C} , \ \ f\mapsto f(C_0 )
$$
\noindent is injective. In this aim, we introduce, for each $\delta
=0,1,2,...$, the set
$$
{\rm Ch}_E^\delta = \{ C\in {\rm Ch}_E \ ; \ d(C,X_F )=\delta\}
$$
\noindent  where $d(C,X_F )$ denotes the combinatorial distance of $C$
to $X_F$:
$$
d(C,X_F )={\rm min}\, \{ d(C,D)\ ; \ D\in X_F\}\ .
$$

\noindent In particular ${\rm Ch}_E^0 = {\rm Ch}_F$. Let $f\in {\mathcal H} (X_F )^{G_F^{\rm
    der}}$. We  then prove that for each $\delta \geqslant 0$,
the restriction of $f$ to ${\rm Ch}_E^{\delta +1}$ depends only on the
restriction of $f$ on ${\rm Ch}_E^\delta$. This follows from  the harmonicity
condition and from a crucial result on the transitivity of the action
of $G_F^{\rm der}$ on the set of chambers of $X_E$ (\cite{BC} Theorem
(5.1)\footnote{The proof of this theorem is due to Fran\c cois
  Court\`es;  see the appendix of \cite{BC}.}). It is now easy to prove by
an inductive argument that the cochain $f$ is known once its values on
${\rm Ch}_E^0 ={\rm Ch}_F$ are known. Since $G_F^{\rm der}$ acts transitively on
${\rm Ch}_F$, $f$ is known once the value $f(C_0 )$ is known and $j$ is
indeed injective.  

\subsection{The tamely ramified case} Conjecture 4.1 was proved by
Fran\c cois Court\`es in the tamely ramified case \cite{Cou2}, i.e. when
$E/F$ is tamely ramified. This case if
much trickier, mainly because,  as we noticed in {\S}1.3, the embedding $X_F\longrightarrow
X_E$ is not simplicial: a chamber of $X_F$ is a union of several
chambers of $X_E$. However the philosophy of Court\`es's approach
remains roughly the same:
\smallskip

 (1) he proves the multiplicity one result by using the model ${\rm
   Hom}_{G_F^{\rm der}}({{\mathbf S}{\mathbf t}}_E ,{\bf 1}) \simeq {\mathcal H} (X_E )^{G_F^{\rm
     der}}$,

 (2) he proves distinction by exhibiting a non zero element of ${\rm
   Hom}_{G_F}({{\mathbf S}{\mathbf t}}_E ,\epsilon_{\rm Prasad})$.
 \smallskip

 For step (1), Court\`es uses an inductive argument similar to that of
 4.4. But a new phenomenon appears : in contrast with the case where
 $E/F$ is unramified the support of a non zero element in ${\mathcal H} (X_E
 )^{G_F^{\rm der}}$ may be quite complicate. In order to analyse this
 support, Court\`es introduces the notion of the {\it anisotropy class} of a
 chamber.
 
 \smallskip

  If $C$ is a chamber of $X_E$ then it belongs to some ${\rm
    Gal}(E/F)$-stable apartment ${\mathcal A}$ of $X_E$ (it is not unique). The
  apartment ${\mathcal A}$ is in turn attached to some ${\rm Gal}(E/F)$-stable maximal $E$-split
torus $T$ of $\mathbb G$. To $T$ one associates its {\it anisotropy class}:
  this is an invariant which describes the ``anitropic part'' of $T$
  as an $F$-torus ($T$ is not necessarily $F$-split). It turns out
  that this anisotropy class does not depend on the choice of $T$;
  this is what Court\`es takes as a definition of the anisotropy class
  of $C$.  Then Court\`es considers two cases.
  \smallskip

   {\it First case: $\mathbb G$ is of type $A_{2n}$.} Write ${\rm Ch}_E^0$
  for the set of chambers of $X_E$ lying in $X_F$. Then any invariant
   non-zero harmonic cochain $f\in {\mathcal H} (X_E )^{\rm G_F^{\rm der}}$ is
   trivial on ${\rm Ch}_E^0$ except on a unique $G_F^{\rm der}$-orbit of
   chambers ${\rm Ch}_c \subset {\rm Ch}_E^0$. Court\`es proves by induction
   that any $f\in  {\mathcal H} (X_E )^{\rm G_F^{\rm der}}$ is entirely
   determined by its restriction to ${\rm Ch}_c$, and multiplicity one
   follows.
   \smallskip

   {\it Second case: $\mathbb G$ is 
     not of type $A_{2n}$, for some integer $n$.} Then any
   $f\in {\mathcal H} (X_E )^{\rm G_F^{\rm der}}$ is trivial on the whole of
   ${\rm Ch}_E^0$ and Court\`es has to find a new starting point for his
   induction argument. It turns out that if $f\in {\mathcal H} (X_E )^{\rm
     G_F^{\rm der}}$ and $C\in {\rm Ch}_E$, the $f(C)=0$ except when $C$
   belongs to a certain anisotropy class of chambers denoted by
   ${\rm Ch}_a$. Court\`es takes as a starting point of his induction the
   set ${\rm Ch}_a^0$ of chambers $C$ of anisotropy class $a$ containing a
   ${\rm Gal}(E/F)$-fixed facet of maximal dimension. He then manages
   to prove that the restriction map
   $$
   {\mathcal H} (X_E )^{\rm  G_F^{\rm der}} \longrightarrow \left\{ f_{\mid {\rm Ch}_a^0}\ ;
   \ f\in {\mathcal H} (X_E )^{ G_F^{\rm der}}\right\}
     $$
    \noindent is injective. He is finally reduced to proving that the space
    of restrictions
    $$
    \left\{ f_{\mid {\rm Ch}_a^0}\ ; \ f\in {\mathcal H} (X_E
    )^{G_F^{\rm der}} \right\}
    $$
    \noindent is one dimensional. This is quite
      technical for the set ${\rm Ch}_a^0$ is not a single $G_F^{\rm
        der}$-orbit in general!

Paul Broussous
\smallskip

paul.broussous{@}math.univ-poitiers.fr
\medskip

D\'epartement de Math\'ematiques

UMR 6086 du CNRS
\smallskip

T\'el\'eport 2 - BP 30179

Boulevard Marie et Pierre Curie

86962 Futuroscope Chasseneuil Cedex

 France

\end{document}